\documentclass[11pt,oneside, final]{amsart}
\copyrightinfo{0}{Iranian Mathematical Society}
\pagespan{1}{\pageref*{LastPage}}
\usepackage{etoolbox,lastpage}
\commby{}
\usepackage{amsmath,amsthm,amscd,amsfonts,amssymb,enumerate}
\usepackage{graphicx}
\usepackage{color}
\usepackage[colorlinks]{hyperref}

\makeatletter \@addtoreset{equation}{section}

\newtheorem{theorem}{Theorem}[section]
\newtheorem{proposition}[theorem]{Proposition}
\newtheorem{lemma}[theorem]{Lemma}

\newtheorem{conjecture}[theorem]{Conjecture}
\newtheorem{question}[theorem]{Question}
\theoremstyle{definition}
\newtheorem{definition}[theorem]{Definition}

\newtheorem{counterexample}[theorem]{Counterexample}
\newtheorem{remark}[theorem]{Remark}
\newtheorem{assumption}[theorem]{Assumption}
\theoremstyle{remark}
\numberwithin{equation}{section}

\begin{document}


\title[Locally constrained flow \& Michael-Simon type
inequalities in hyperbolic space]{Locally constrained flows and
Michael-Simon type inequalities in hyperbolic space
$\mathbb{H}^{n+1}$}

\author[J. Cui]{Jingshi Cui}

\author[P. Zhao]{Peibiao Zhao$^*$}

\thanks{$^*$Corresponding author}
%
\begin{abstract}
Brendle \cite{B21} successfully establishes the sharp Michael-Simon inequality for mean curvature on Riemannian manifolds with nonnegative sectional curvature ($\mathcal{K} \geq 0$), and the proof relies on the Alexandrov-Bakelman-Pucci method. Nevertheless, this result cannot be extended to hyperbolic space $\mathbb{H}^{n+1}$ ($\mathcal{K} = -1$), as demonstrated by Counterexample \ref{ex1.7}.

In the present paper, we propose Conjectures \ref{con1.8} and {\ref{con1.9}} concerning the hyperbolic version of the sharp Michael-Simon type inequality for $k$-th mean curvatures. However, the proof method in \cite{B21} failed to verify the validity of these conjectures. Recently, the authors \cite{CZ24} proved Conjectures \ref{con1.8} and {\ref{con1.9}} only for $h$-convex hypersurfaces by means of the Brendle-Guan-Li's flow. 

This paper aims to utilize other types of curvature flows to prove Conjectures \ref{con1.8} and {\ref{con1.9}} for hypersurfaces with weaker convexity conditions. For $k = 1$, we first investigate a new locally constrained mean curvature flow (\ref{1.9}) in $\mathbb{H}^{n+1}$ and prove its longtime existence and exponential convergence. Then, the sharp Michael-Simon type inequality for mean curvature of starshaped hypersurfaces in $\mathbb{H}^{n+1}$ is confirmed through the flow (\ref{1.9}). For $k \geq 2$, the sharp Michael-Simon inequality for $k$-th mean curvatures of starshaped, strictly $k$-convex hypersurfaces in $\mathbb{H}^{n+1}$ is proven using the locally constrained inverse curvature flow (\ref{1.12}) introduced by Scheuer and Xia \cite{SX19}.
\end{abstract}

\maketitle
%

\section{ Introduction}

\noindent  The Michael-Simon inequality on the generalized submanifolds immersed into Euclidean space $\mathbb{R}^{N}$ was first formulated by Michael and Simon \cite{MS1973}. The specific form of the inequality is as follows
\begin{theorem}{(\cite{MS1973})}\label{th1.1}
    Let $i: M^{n} \rightarrow \mathbb{R}^{N}$ be an isometric immersion $(N>n)$. Let $U$ be an open subset of $M^n$. For a function $\phi   \in C_{c}^{\infty}(U)$, there exists a constant $C$, such that
    \begin{align}\label{1.1}
        \int_{M^{n}}(|H| \cdot|\phi|+|\nabla \phi|) d \mu_{M^{n}} \geq C\left(\int_{M^{n}}|\phi|^{\frac{n}{n-1}} d \mu_{M^{n}}\right)^{\frac{n-1}{n}},
    \end{align}
    where $H$ is the mean curvature of $M^n$.
\end{theorem}
There are many important applications of the Michael-Simon type inequality. First, this type inequality is an effective tool for obtaining a priori estimates of solutions to elliptic partial differential equations and for investigating the regularity of these solutions. For example, X. Cabr$\acute{\rm e} $ \cite{Ca10} used this type inequality to prove the regularity theory for semilinear elliptic equations. Second, since the Michael-Simon type inequality involves integrals of mean curvature, it allows for solving problems related to hypersurfaces with prescribed mean curvature. For instance, this type inequality can be applied to study the regularity of surfaces with prescribed mean curvature \cite{BG73, DHT10}. Finally, the Michael-Simon type inequality is also valuable in the study of geometric flows. For instance, Evans and Spruck \cite{ES92} used this type inequality to establish an analog of Brakke’s ``clearing out'' lemma for generalized mean curvature flow. In summary, the Michael-Simon type inequality plays a crucial role not only in the study of elliptic equations but also in geometric analysis and the investigation of geometric flows. Therefore, it is necessary to explore the Michael-Simon type inequality for the $k$-th mean curvature, as well as the sharp form of this type inequality. Moreover, it is also important to investigate the Michael-Simon type inequality in Riemannian manifolds.

It is widely known that Sun-Yung Alice Chang and Yi Wang \cite{CW2013}, utilized the optimal transport mapping method to successfully establish the Michael-Simon type inequality for $k$-th mean curvatures. Subsequently, in 2021, Brendle \cite{B19} established a sharp form of the Michael-Simon type inequality, where the best constant $C$ corresponds to the isometric constant in $\mathbb{R}^{n+1}$.
\begin{theorem}{(\cite{B19})}\label{th1.2}
    Let $M$ be a compact hypersurface in $\mathbb{R}^{n+1}$ (possibly with boundary $\partial M$), and let $f$ be a positive smooth function on $M$. Then
    \begin{align}\label{1.2}
       \int_{M} \sqrt{\left|\nabla^{M} f\right|^{2}+f^{2} H^{2}}+\int_{\partial M} f \geq n\omega_{n}^{\frac{1}{n}}\left(\int_{M} f^{\frac{n}{n-1}}\right)^{\frac{n-1}{n}},
    \end{align}
    where $H$ is the mean curvature of $M$ and $\omega_{n}$ is the area of the unit sphere $\mathbb{S}^{n}$ in $\mathbb{R}^{n+1}$.
    Moreover, if the equality holds, then $f$ is constant and $M$ is a flat disk.
\end{theorem}

\begin{remark}
    (1) If $M$ is a compact minimal hypersurface in $\mathbb{R}^{n+1}$, then by the Michael-Simon type inequality (\ref{1.2}), $M$ satisfies the sharp isoperimetric inequality
    \begin{align*}
        |\partial M| \geq n\omega_{n}^{\frac{1}{n}}|M|^{\frac{n-1}{n}}.
    \end{align*}
    (2) If $M$ is a closed hypersurface (without boundary) and $f$ be of constant, then the inequality (\ref{1.2}) gives the famous Minkowski inequality in $\mathbb{R}^{n+1}$, i.e.,
    \begin{align*}
        \int_{M} |H| d\mu \geq n \omega^{\frac{1}{n}}_{n}|M|^{\frac{n-1}{n}}.
    \end{align*}
\end{remark}

Moreover, in 2022, Brendle \cite{B21} established the Michael-Simon type inequality for mean curvature in Riemannian manifolds with nonnegative sectional curvature. 

\begin{theorem}{(\cite{B21})\label{th1.5}}
    Let $N$ be a complete noncompact manifold of dimension $n+m$ with nonnegative sectional curvatures. Let $M$ be a compact submanifold of $N$ of dimension $n$ (possibly with boundary $\partial M$), and let $f$ be a positive smooth function on $M$. If $m\geq 2 $, then
    \begin{align}\label{1.3}
        \int_{M} \sqrt{\left|\nabla^{\Sigma} f\right|^{2}+f^{2} |H|^{2}}+\int_{\partial M} f \geq n\left(\frac{(n+m)\omega_{n+m}}{m\omega_{m}}\right)^{\frac{1}{n}}\theta^{\frac{1}{n}}\left(\int_{M} f^{\frac{n}{n-1}}\right)^{\frac{n-1}{n}},
    \end{align}
    where $\theta$ denotes the asymptotic volume ratio of $N$ and $H$ denotes the mean curvature vector of $M$.
\end{theorem}

\begin{remark}\label{rem1.6}
    (1) When $m=2$, we have $(n+2) \omega_{n+2}=2 \omega_{2} \cdot \omega_{n}$. Then, the inequality (\ref{1.3}) is as follows
    \begin{align}\label{1.4}
        \int_{M} \sqrt{\left|\nabla^{\Sigma} f\right|^{2}+f^{2} |H|^{2}}+\int_{\partial M} f \geq n\omega_{n}^{\frac{1}{n}}\theta^{\frac{1}{n}}\left(\int_{M} f^{\frac{n}{n-1}}\right)^{\frac{n-1}{n}}.
    \end{align}

    (2) The inequality (\ref{1.4}) also holds for the setting of codimension one. Since $M$ can be regarded as a submanifold of the $(n+2)$-dimensional manifold $N\times\mathbb{R}$, the asymptotic volume ratios of $N\times\mathbb{R}$ and $N$ are equal.

    (3) Note that the inequalities (\ref{1.3}) and (\ref{1.4}) depend on the asymptotic volume ratio $\theta$ of the ambient manifold $N$. From the Bishop-Gromov volume comparison theorem \cite{BC64,MG81}, it follows that the following two conclusions hold: (i) $\theta \leq 1$ when $(N,g_{N})$ is a complete Riemannian manifold with Ricci curvature bounded below; (ii) $\theta > 1$ when the sectional curvature of $(N,g_{N})$ is negative.
\end{remark}

Given the significant importance of the Michael-Simon type inequality in the relevant mathematical fields, it is both necessary and urgent to conduct an study of the Michael-Simon type inequality on Riemannian manifolds with negative sectional curvature. As everyone knows, the hyperbolic space $\mathbb{H}^{n+1}$ is a simply connected, complete Riemannian ($n+1$)-manifold with $\mathcal{K}=-1<0$. Therefore, we first attempt to establish the Michael-Simon type inequality in $\mathbb{H}^{n+1}$. A natural question arises
\begin{question}\label{qu1.4}
    Does the Michael-Simon type inequality in \cite{B21} also hold in $\mathbb{H}^{n+1}$ ?
\end{question}

We find that when $m=1$ and $f=const$, there exist counterexamples to the inequality (\ref{1.4}) in hyperbolic space $\mathbb{H}^{n+1}$.
\begin{counterexample}\label{ex1.7}
    Let $M$ be a geodesic sphere $B_{r}$ in $\mathbb{H}^{n+1}$ where $r\in (0,+\infty)$. Assume that the hyperbolic space $\mathbb{H}^{n+1}$ is equipped with the measure $d\mu_{\mathbb{H}^{n+1}}=\lambda^{n}drd\mu_{\mathbb{S}^{n}}$, then $H\big|_{B_{r}} = \frac{n\lambda^{'}}{\lambda}$ and the inequality (\ref{1.4}) becomes
    \begin{align*}
        \int_{B_{r}} \frac{\lambda^{'}}{\lambda} d\mu_{\mathbb{H}^{n+1}} \geq \theta^{\frac{1}{n}} \omega_{n}^{\frac{1}{n}}|B_{r}|^{\frac{n-1}{n}}.
    \end{align*}
    Simplifying the above inequality, we have
    \begin{align*}
        \theta \leq \int^{r}_{0} (\lambda^{'})^{n} ds = \int^{r}_{0} \left(\frac{e^{s}+e^{-s}}{2}\right)^{n} ds < \int^{r}_{0} e^{ns} ds =\frac{e^{nr}-1}{n}.
    \end{align*}
    According to the Bishop-Gromov volume comparison theorem, in the hyperbolic space $\mathbb{H}^{n+1}$, $\theta >1$. However, if $r\in (0,\frac{ln(n+1)}{n} ]$, then $\theta < 1$, which is clearly in contradiction with the Bishop-Gromov volume comparison theorem.
\end{counterexample}
Counterexample \ref{ex1.7} provides an answer to Question \ref{qu1.4}, namely, the conclusion in \cite{B21} does not hold in hyperbolic space $\mathbb{H}^{n+1}$. Therefore, we make conjectures for the Michael-Simon type inequality in $\mathbb{H}^{n+1}$.
\begin{conjecture}\label{con1.8}
    Let $M$ be a compact hypersurface in $\mathbb{H}^{n+1}$ (possibly with boundary $\partial M$) and let $f$ be a positive smooth function on $M$. Then,
    \begin{align}\label{1.5}
        \int_{M} \lambda^{'} \sqrt{f^{2} E_{1}^{2}+|\nabla^{M} f|^{2}} -\int_{M}\left\langle \bar{\nabla}\left(f\lambda^{'}\right),\nu \right\rangle  +\int_{\partial M}  f
        \geq
        C\left(\int_{M}f^{\frac{n}{n-1}}\right)^{\frac{n-1}{n}} .
    \end{align}
    In particular, when $M$ is closed and $f=const$, there holds
    \begin{align}\label{1.6}
        \int_{M}(\lambda^{'}|E_{1}|-u)  \geq C|M|^{\frac{n-1}{n}},
    \end{align}
    where $E_{1}:=E_{1}(\kappa)$ and $\nu$ are the (normalized) $1$-th mean curvature and the unit outward normal of $M$ respectively, $u$ and $|M|$ are the support function and the area of $M$ respectively, and $\bar{\nabla}$ is the Levi-Civita connection with respect to the metric $\bar{g}$ on $\mathbb{H}^{n+1}$.  
\end{conjecture}
Note that the inequality (\ref{1.6}) is precisely the Minkowski type inequality established by Brendle, Hung, and Wang in \cite{BHW16} for mean convex and starshaped hypersurfaces. Similarly, we propose a conjecture for the Michael-Simon type inequality for $k$-th mean curvatures in hyperbolic space $\mathbb{H}^{n+1}$.
\begin{conjecture}\label{con1.9}
    Let $M$ be a compact hypersurface in $\mathbb{H}^{n+1}$ (possibly with boundary $\partial M$) and let $f$ be a positive smooth function on $M$. For $1\leq k\leq n$, there holds
    \begin{align}\label{1.7}
        \int_{M} \lambda^{'} \sqrt{f^{2}E_{k}^{2} +|\nabla^{M} f|^{2} E_{k-1}^{2}} -\int_{M}\left\langle \bar{\nabla}\left(f\lambda^{'}\right),\nu \right\rangle \cdot |E_{k-1}| &+\int_{\partial M}  f \cdot|E_{k-1}| \notag \\
        &\geq
        C\left(\int_{M} f^{\frac{n-k+1}{n-k}}\cdot |E_{k-1}|\right)^{\frac{n-k}{n-k+1}}. 
    \end{align}
In particular, when $M$ is closed and $f=const$, there holds
    \begin{align}\label{1.8}
        \int_{M}(\lambda^{'}|E_{k}|-u|E_{k-1}|) \geq C \left(\int_{M} |E_{k-1}|\right)^{\frac{n-k+1}{n-k}},
    \end{align}
    where $E_{k}:=E_{k}(\kappa)$ is the (normalized) $k$-th mean curvatures of $M$.
\end{conjecture}
When $M$ is a closed, $k$-convex hypersurface with $k\geq 2$, applying the Minkowski equality (\ref{2.9}) in $\mathbb{H}^{n+1}$, the inequality (\ref{1.8}) implies
\begin{align*}
    \int_{M}(\lambda^{'}E_{k}-\lambda^{'}E_{k-2}) \geq C \left(\int_{M} E_{k-1}\right)^{\frac{n-k+1}{n-k}},
\end{align*}
which represents a geometric inequality for weighted curvature integrals $W^{\lambda^{'}}_{k+1} = \int_{M}\lambda^{'}E_{k}$. The geometric inequality related to $W^{\lambda^{'}}_{k+1}$ can be found in \cite{HLW20}.

Nowadays, the study of curvature flows is in full swing, and has achieved extremely fruitful research results. Notably, the curvature flow approach is a powerful tool for proving new geometric inequalities, especially constrained curvature flows. As we know that the research on constrained curvature flows can be divided into two categories: globally constrained flows and locally constrained flows. For globally constrained flows, Wang and Xia \cite{WangXia14} utilized the volume-preserving mean curvature flow, which was initially proposed by Cabezas-Rivas and Miquel \cite{CM07}, developed the quermassintegral inequalities for $h$-convex hypersurfaces and solved the isoperimetric problem in $\mathbb{H}^{n+1}$. See also \cite{BP17,GLW17,Mak12} for studies on the globally constrained flow in $\mathbb{H}^{n+1}$. Regarding locally constrained flows, Hu, Li, and Wei \cite{HLW20} not only provided new proofs of the quermassintegral inequalities but also derived new geometric inequalities in $\mathbb{H}^{n+1}$. Moreover, Hu and Li \cite{HL21} proved the weighted geometric inequality in hyperbolic space $\mathbb{H}^{n+1}$ by virtue of the locally constrained weighted volume preserving flows. Other results on locally constrained curvature flows, see \cite{JX19, LS20, SWX18, WX19, WX, GL2015, GL19}.

Although the author's recent work \cite{CZ24} uses the Brendle-Guan-Li's flow to prove Conjectures \ref{con1.8} and \ref{con1.9}, the result requires the hypersurfaces to be $h$-convex, i.e., $\kappa_{i} \geq 1$. If the convexity condition of the hypersurfaces can be weakened, the scope of application of the inequality can be broadened. This paper aims to prove Conjectures \ref{con1.8} and \ref{con1.9} using other types of curvature flows, thereby enabling hypersurfaces with weaker convexity conditions to satisfy the inequality.

In the first part of the paper, we introduce a new locally constrained mean curvature flow for starshaped hypersurfaces in $\mathbb{H}^{n+1}$, and apply it to establish the Michael-Simon type inequality for mean curvature, i.e., to prove Conjecture \ref{con1.8}.

Let $X:\Sigma \times [0,T) \to \mathbb{H}^{n+1}$ be a smooth family of embedding such that $\Sigma_{t}=X(\Sigma,t)$ are smooth, closed starshaped hypersurfaces in $\mathbb{H}^{n+1}$. The new locally constrained mean curvature flow is 
\begin{equation}\label{1.9}
    \begin{cases}
        \frac{\partial }{\partial t}X(x,t)=-\frac{f}{v}H \nu(x,t) - \frac{n}{n-1}\frac{1}{v}\frac{\partial f}{\partial X},\\
        X(\cdot,0)=X_{0}(\cdot),
     \end{cases}
\end{equation}
where $f \in C^{\infty}(\Sigma_{t})$ is a positive function, $\nu(x,t)$ and $H$ are the unit outer normal and the mean curvature of $\Sigma_{t}$ respectively, and $v=\sqrt{1+\lambda^{-2}|Dr|^{2}}$, with $D$ denotes Levi-Civita connection on $\mathbb{S}^{n}$.

It is known that starshaped hypersurfaces $\Sigma_{t}$ can be parametrized by the positive radial function $r(\cdot,t): \mathbb{S}^{n}\to \mathbb{R}^{+}$ on $\mathbb{S}^{n}$, i.e.,
\begin{align*}
    \Sigma_{t}=X(\cdot,t)=\left\{\left(r(\xi,t),\xi\right)\in \mathbb{R}^{+}\times\mathbb{S}^{n}|\xi \in\mathbb{S}^{n}\right\},
\end{align*}
and $\Sigma_{0}=X(\cdot,0)=\left\{\left(r_{0}(\xi),\xi\right)\in \mathbb{R}^{+}\times\mathbb{S}^{n}|\xi \in\mathbb{S}^{n}\right\}$, with $r_{0}(\xi)=r(\xi,0)$. Since $f \in C^{\infty}(\Sigma_{t})$, it follows that there exists a function $\phi: \mathbb{S}^{n} \to \mathbb{R}^{+} $ such that $f(X) = f(r(\xi,t),\xi)= \phi (\xi)$.

Next we prove that the longtime existence and convergence of the locally constrained mean curvature flow (\ref{1.9}), provided that $f$ satisfies the following assumption.
\begin{assumption}\label{as1.10}
    Let $f(X) = \Phi_{1} \circ r(\xi,t) $ and $\widehat{\Phi}_{1}(r):=n\Phi_{1}\frac{\lambda^{'}}{\lambda^{2}}+\frac{n}{n-1} \frac{\partial \Phi_{1}}{\partial r}\frac{1}{\lambda}$. The function $\widehat{\Phi}_{1}(r)$ is monotonically increasing with respect to $r$, where $r>0$, and there exists a zero point for $\widehat{\Phi}_{1}(r)$.
\end{assumption}

\begin{remark}
    It is not difficult to verify that there exists $\Phi_{1}(r)$ and $\widehat{\Phi}_{1}(r)$ satisfies Assumption \ref{as1.10}. For example, let
    \begin{align*}
        \Phi_{1}(r) = \lambda^{-n+1}\left[\int\frac{n-1}{n}(e^{-1} - e^{-r})\lambda^{n} dr +C\right]
    \end{align*}
    where the constant $C$ is chosen as follows: for $r >1$, let $C=0$; for $0< \epsilon \leq r \leq1 $, let $C \geq  (1-\epsilon) \frac{n-1}{n}\frac{e-1}{e}\lambda^{n}(1)$. With this choice of  $C$, $\Phi_{1}(r)$ is positive for all $r>0$. Moreover, we have $\widehat{\Phi}_{1}(r) = e^{-1} - e^{-r}$, $\frac{\partial \widehat{\Phi}_{1}}{\partial r} = e^{-r} >0$ and $\widehat{\Phi}_{1}(1) = 0$.
\end{remark}
\begin{theorem}\label{th1.11}
    Let $X_{0}:\Sigma \to \mathbb{H}^{n+1}(n\geq 2)$ be a smooth embedding of closed hypersurface $\Sigma$ in $\mathbb{H}^{n+1}$ such that $\Sigma_{0}=X_{0}(\Sigma)$ is starshaped with respect to the origin, and assume that $f$ satisfies Assumption \ref{as1.10}. Then the flow (\ref{1.9}) has a unique smooth solution $\Sigma_{t}=X(\Sigma,t)$ for all time $t\in [0,+\infty)$. Moreover, $\Sigma_{t}$ converges exponentially to a geodesic sphere $B_{r_{\infty}}$ centered at the origin as $t\to +\infty $ in the $C^{\infty}$-topology.
\end{theorem}

Let $M_{0}= \Sigma_{0}$ ($M_{0}=\Sigma_{0} \cup \partial \Sigma_{0}$ ) be a smooth starshaped and compact hypersurface (possibly with boundary). For any $ X^{*} \in M_{0}$, we have  $X^{*}=(r^{*}(\xi), \xi)$, $r^{*}$ is the distance function of $M_{0}$ and  
\begin{equation*}
    r^{*}(\xi) =
    \begin{cases}
        r_{0}(\xi),\qquad x \in \Sigma_{0},
       \\
       \bar{r}(\xi), \qquad x \in \partial \Sigma_{0}.
    \end{cases}
\end{equation*}
If $f\in C^{\infty}(M_{0})$, then $f(X^{*})= f(\left(r^{*}(\xi),\xi\right))$ and $f(X^{*})\big|_{\Sigma_{0}} = f(X(\cdot,0)) = f (r_{0}(\xi),\xi)$. By using the locally constrained mean curvature flow (\ref{1.9}), the sharp Michael-Simon type inequality for mean curvature in $\mathbb{H}^{n+1}$  is established, and the best constant in (\ref{1.5}) is determined.
\begin{theorem}\label{th1.12}
    Let $M_{0}$ be a smooth compact hypersurface in $\mathbb{H}^{n+1}$ (possibly with boundary $\partial M_{0}$) that is starshaped with respect to the origin, and let $f \in C^{\infty}(M_{0})$ be a positive function of the form $f=\Phi_{1} \circ r$ satisfying Assumption \ref{as1.10}. Then, there holds
    \begin{align}\label{1.10}
        \int_{M_{0}} \lambda^{'} \sqrt{f^{2} E_{1}^{2}+|\nabla^{M} f|^{2}} -\int_{M_{0}}\left\langle \bar{\nabla}\left(f\lambda^{'}\right),\nu \right\rangle  +\int_{\partial M_{0}}  f
        \geq
        \omega_{n}^{\frac{1}{n}}\left(\int_{M_{0}}f^{\frac{n}{n-1}}\right)^{\frac{n-1}{n}}.
    \end{align}
    In particular, if $M_{0}$ is closed, we have
    \begin{align}\label{1.11}
        \int_{M_{0}} \lambda^{'} \sqrt{f^{2} E_{1}^{2}+|\nabla^{M} f|^{2}} -\int_{M_{0}}\left\langle \bar{\nabla}\left(f\lambda^{'}\right),\nu \right\rangle 
        \geq
        \omega_{n}^{\frac{1}{n}}\left(\int_{M_{0}}f^{\frac{n}{n-1}}\right)^{\frac{n-1}{n}}. \tag{1.10 ${'}$}
    \end{align}
    Equality holds in (\ref{1.11}) if and only if $M_{0}$ is a geodesic sphere centered at the origin.
\end{theorem}

In the second part of this paper, we use the locally constrained inverse curvature type flow to establish the Michael-Simon type inequality for $k$-th mean curvatures, i.e., to prove Conjecture \ref{con1.9}.

Let $X_{0}:\Sigma \to \mathbb{H}^{n+1}$ be a smooth embedding such that $\Sigma_{0}$ is a closed, starshaped, and strictly $k$-convex hypersurface in hyperbolic space $\mathbb{H}^{n+1}$. In \cite{SX19}, Scheuer and Xia introduced a family of smooth embeddings $X:\Sigma \times [0,T) \to \mathbb{H}^{n+1}$ that satisfy
\begin{equation}\label{1.12}
    \begin{cases}
        \frac{\partial}{\partial t}X(x,t)=\left(\frac{E_{k-2}(\kappa)}{E_{k-1}(\kappa)}-\frac{u}{\lambda^{'}}\right)\nu(x,t),\quad k=2,\cdots ,n,\\
        X(\cdot,0)=X_{0}(\cdot),
    \end{cases}
\end{equation}
where $\nu(x,t)$ and $\kappa=(\kappa_{1},\cdots,\kappa_{n})$ are the unit outer normal and the principal curvatures of $\Sigma_{t}=X(\Sigma,t)$ respectively. They proved the following convergence result of the flow (\ref{1.12}).
\begin{theorem}\label{th1.13}(\cite{SX19})
    Let $X_{0}(\Sigma)$ be the smooth embedding of a closed n-dimensional manifold $\Sigma$ in $\mathbb{H}^{n+1}$, such that $\Sigma_{0}=X_{0}(\Sigma)$ is starshaped and strictly $k$-convex along $X_{0}(\Sigma)$. Then any solution $\Sigma_{t}=X_{t}(\Sigma)$ of the flow (\ref{1.12}) exists for $t\in[0,+\infty)$. Moreover, $\Sigma_{t}$ is starshaped and strictly  $k$-convex for each $t> 0$ and it converges to a geodesic sphere $B_{R}$ centered at the origin in the $C^{\infty}$-topology as $t \to +\infty$.
\end{theorem}

Moreover, as an application of the flow (\ref{1.12}), we obtain the sharp Michael-Simon type inequalities for $k$-th mean curvatures, where $k\geq 2$. Before presenting the detailed results, we make the following assumptions about $f$.
\begin{assumption}\label{as1.14}
    Let $f(X) = \Phi_{2} \circ r(\xi,t)$ and $\widehat{\Phi}_{2}(r):= \Phi_{2}^{\frac{n-k+1}{n-k}}(r) $.  
    \\
    (1) $\widehat{\Phi}_{2}$ is the solution of the following second order non-homogeneous differential equation
    \begin{align}\label{1.13} 
        \frac{1}{n}\Delta_{\Sigma_{t}} \widehat{\Phi}_{2} - \widehat{\Phi}_{2} - \frac{1}{k-1}\frac{\lambda^{'}}{\lambda}\frac{\partial \widehat{\Phi}_{2}}{\partial r}\left(1+\lambda^{-2}|Dr|^{2}\right) =- \frac{m}{k-1}\frac{E_{k-1}}{E_{k-2}}\lambda (\lambda^{'})^{m-k}, 
    \end{align}
    where $k+1 \leq m \leq n$.
    \\
    (2) $(\lambda^{'})^{-1}\widehat{\Phi}_{2}$ is monotonically increasing with respect to $\lambda^{'}$.
\end{assumption}
\begin{remark}
    By the standard theory of ordinary differential equations, the solution to the equation (\ref{1.13}) exists. It can be shown that the set of elements satisfying Assumption \ref{as1.14} is non-empty. In particular, when $\Sigma_{t}$ is a geodesic sphere $B_{R}$, $\widehat{\Phi}_{2} = (\lambda^{'})^{m-k+1}$ satisfies Assumption \ref{as1.14}.
\end{remark}

\begin{theorem}\label{th1.15}
    Let $M_{0}$ be a smooth, compact, starshaped and strictly $k$-convex hypersurface in $\mathbb{H}^{n+1}$ (possibly with boundary $\partial M_{0}$), and $\Omega_{0}$ be the domain enclosed by $M_{0}$. Assume that $f$ satisfies Assumption \ref{as1.14}. Then, for any $2\leq k \leq n$, there holds
    \begin{align}\label{1.14}
        \int_{M_{0}} \lambda^{'} \sqrt{f^{2}E_{k}^{2} +|\nabla^{M_{0}} f|^{2} E_{k-1}^{2}} &-\int_{M_{0}}\left\langle \bar{\nabla}\left(f\lambda^{'}\right),\nu \right\rangle \cdot E_{k-1} +\int_{\partial M_{0}}  f \cdot E_{k-1} \notag \\
        &\geq
        \left(p_{k}\circ h_{0}^{-1}\left(W^{\lambda^{'}}_{0}(\Omega_{0})\right) \right)^{\frac{1}{n-k+1}}\left(\int_{M_{0}}f^{\frac{n-k+1}{n-k}}\cdot E_{k-1}\right)^{\frac{n-k}{n-k+1}}.
    \end{align}
    In particular, if $M_{0}$ is closed, we have
    \begin{align}\label{1.15}
        \int_{M_{0}} \lambda^{'} \sqrt{f^{2}E_{k}^{2} +|\nabla^{M_{0}} f|^{2} E_{k-1}^{2}} &-\int_{M_{0}}\left\langle \bar{\nabla}\left(f\lambda^{'}\right),\nu \right\rangle \cdot E_{k-1} \notag \\
        &\geq
        \left(p_{k}\circ h_{0}^{-1}\left(W^{\lambda^{'}}_{0}(\Omega_{0})\right) \right)^{\frac{1}{n-k+1}}\left(\int_{M_{0}}f^{\frac{n-k+1}{n-k}}\cdot E_{k-1}\right)^{\frac{n-k}{n-k+1}}, \tag{1.13 ${'}$}
    \end{align}
    where $p_{k}(R)=\omega_{n}f^{\frac{n+1-k}{n-k}}(\lambda^{'})^{k-1}\lambda^{n-k+1}(R)$, $W^{\lambda^{'}}_{0}(\Omega_{0})= (n+1)\int_{\Omega_{0}} \lambda^{'} dvol$,  $h_{0}(R)=W^{\lambda^{'}}_{0}(B^{n+1}_{R})$ and $h_{0}^{-1}$ is the inverse function of $h_{0}$. Equality holds in (\ref{1.15}) if and only if $M_{0}$ is a geodesic sphere centered at the origin.
\end{theorem}

The paper is organized as follows. In Section 2, we review the geometry of starshaped hypersurfaces in hyperbolic space $\mathbb{H}^{n+1}$, present some basic properties of normalized elementary symmetric functions, and derive the corresponding evolution equations. In Section 3, we establish a priori estimates for the locally constrained mean curvature flow (\ref{1.9}), prove its longtime existence and convergence. The sharp Michael-Simon type inequalities in Theorem \ref{th1.12} and Theorem \ref{th1.15} are proved in Section 4.
\section{ Preliminaries}

\noindent In this section, we first review starshaped hypersurfaces in hyperbolic space $\mathbb{H}^{n+1}$, then state some properties of normalized elementary symmetry functions, and finally derive the evolution equations along the general flow (\ref{2.11}).

\subsection{ Starshaped hypersurfaces in a hyperbolic space.}
\
\vglue-10pt
 \indent

 In this paper, the hyperbolic space $\mathbb{H}^{n+1}$ is regarded as a warped product manifold $\mathbb{R}^{+} \times \mathbb{S}^{n}$ equipped with the metric
\begin{align*}
    \bar{g} =dr^{2}+\lambda(r)^{2}\sigma,
\end{align*}
where $\sigma$ is the standard metric on the unit sphere $\mathbb{S}^{n}\subset \mathbb{R}^{n+1}$. Let $\Sigma$ be a closed smooth hypersurface in $\mathbb{H}^{n+1}$. Define $g_{ij}$, $h_{ij}$ and $\nu$ as the induced metric, the second fundamental form and the unit outward normal vector of $\Sigma$, respectively. The principal curvatures $\kappa=(\kappa_{1},\cdots,\kappa_{n})$ of $\Sigma$ are given by the eigenvalues of the Weingarten matrix $\mathcal{W}=(h^{i}_{j})=(g^{ik}h_{kj})$, where $g^{ik}$ is the inverse matrix of $g_{ik}$.

\begin{lemma}(\cite{GL2015})\label{lem2.1}
    Let $(\Sigma,g)$ be a smooth hypersurface in $\mathbb{H}^{n+1}$. Denote
    \begin{align*}
        \Gamma (r)=\int^{r}_{0} \lambda(y)dy=\lambda^{'}(r)-1.  
    \end{align*}
    Then, $\Gamma|_{\Sigma}$ satisfies
    \begin{align}
        \nabla_{i}\Gamma &=\nabla_{i} \lambda^{'} =\left\langle \lambda \partial r, e_{i}\right\rangle, \label{2.1} \\
        \nabla_{i}\nabla_{j}\Gamma &= \nabla_{i}\nabla_{j}\lambda^{'}=\lambda^{'}g_{ij}-uh_{ij}, \label{2.2}
    \end{align}
    where $\left\{e_{1},\cdots,e_{n}\right\}$ is a basis of the tangent space of $\Sigma$.
\end{lemma}

As $\Sigma \subset \mathbb{H}^{n+1}$ is a smooth, closed starshaped hypersurface with respect to the origin, then the support function $u>0$ everywhere on $\Sigma$. Moreover, $\Sigma$ can be represented as a radial graph, i.e., $\Sigma=\{ (r(\xi),\xi), \xi \in \mathbb{S}^{n}\}$. Let $\xi=(\xi^{1},\cdots,\xi^{n})$ be a local coordinate of $\mathbb{S}^{n}$, $\partial_{i}=\partial_{{\xi}^{i}}$, $r_{i}=D_{i}r$, with $D$ is the Levi-Civita connection on $(\mathbb{S}^{n},\sigma)$. We introduce a new function $\varphi:\mathbb{S}^{n} \to \mathbb{R}$ by
\begin{align*}
    \varphi(\xi)=\chi(r(\xi)),
\end{align*}
where $\chi$ is a positive smooth function that satisfies $\frac{\partial }{\partial r}\chi=\frac{1}{\lambda(r)}$. Thus
\begin{align*}
    \varphi_{i}:=D_{i}\varphi=\frac{r_{i}}{\lambda(r)}.
\end{align*}
The induced metric $g_{ij}$ and its inverse matrix $g^{ij}$ have the following forms
\begin{align*}
    &g_{i j}=\lambda^{2}\sigma_{i j}+r_{i} r_{j}= \lambda^{2}\left(\sigma_{i j}+ \varphi_{i} \varphi_{j}\right)  ,\\
    &g^{i j}=\frac{1}{\lambda^{2}}\left(\sigma^{ij}-\frac{r^{i} r^{j}}{\lambda^{2}+|Dr|^{2}}\right)=\frac{1}{\lambda^{2}}\left(\sigma^{ij}-\frac{\varphi^{i} \varphi^{j}}{v^{2}}\right),
\end{align*}
where $r^{i}=\sigma^{ik}r_{k}$, $\varphi^{i}=\sigma^{ik}\varphi_{k}$ and $v=\sqrt{1+\lambda^{-2}|Dr|^{2}}=\sqrt{1+|D\varphi|^{2}}$.
The unit outer normal $\nu$ and the support function $u$ can be expressed by
\begin{align*}
    \nu= \frac{1}{v}\left(\partial_{r} - \frac{r_{i}}{\lambda^{2}}\partial_{i}\right) =\frac{1}{v}\left(\partial_{r} - \frac{\varphi_{i}}{\lambda}\partial_{i}\right), \qquad u=\left\langle \lambda\partial_{r}, \nu \right\rangle =\frac{\lambda}{v}.
\end{align*}
The second fundamental form $h_{ij}$, the Weingarten matrix $h^{i}_{j}$ and the mean curvature of $\Sigma$ are as follows (see e.g., \cite{G11})
\begin{align}
    &h_{ij}=\frac{\lambda^{'}}{\lambda v}g_{ij}-\frac{\lambda}{v}\varphi_{ij} ,\label{2.3}\\
    &h^{i}_{j}=g^{ik}h_{kj}= \frac{\lambda^{'}}{\lambda v}\delta^{i}_{j}-\frac{1}{\lambda v}\left(\sigma^{ik}-\frac{\varphi^{i}\varphi^{k}}{v^{2}}\right)\varphi_{kj}, \label{2.4}\\
    &H=\frac{n\lambda^{'}}{\lambda v}-\frac{1}{\lambda v}\left(\sigma^{ik}-\frac{\varphi^{i}\varphi^{k}}{v^{2}}\right)\varphi_{ki}.\label{2.5}
\end{align}

\subsection{ Normalized elementary symmetric
functions }
\
\vglue-10pt
 \indent

For each $l=1, \cdots, n$, the normalized $l$-th elementary symmetric functions for $\kappa=\left(\kappa_{1}, \cdots,
\kappa_{n}\right)$ are
\begin{align*}
    E_{l}(\kappa)=\binom{n}{l}^{-1} \sigma_{l}(\kappa)=\binom{n}{l}^{-1} \sum_{1 \leq i_{1}<\ldots<i_{l} \leq n} \kappa_{i_{1}} \cdots \kappa_{i_{l}},
\end{align*}
and we can set $E_{0}(\kappa)=1$ and $E_{l}(\kappa)=0$ for $l>n$. If $A=[A_{ij}] $ is an $n \times n$ symmetric matrix and $\kappa=\kappa(A)=(\kappa_{1},\cdots,\kappa_{n})$ are the eigenvalues of $A$, then $E_{l}(A)=E_{l}(\kappa(A))$ can be expressed as
\begin{align*}
    E_{l}(A)=\frac{(n-l)!}{n!} \delta_{i_{1} \ldots i_{l}}^{j_{1} \ldots j_{l}} A_{i_{1} j_{1}} \cdots A_{i_{l} j_{l}}, \quad l=1, \ldots, n,
\end{align*}
where $\delta_{i_{1} \ldots i_{l}}^{j_{1} \ldots j_{l}}$ is a generalized Kronecker delta. See \cite{G2013} for details. We now recollect some properties of normalized $l$-th elementary symmetric functions.
\begin{lemma}(\cite{G2013})\label{lem2.2}
    Let $\dot{E}_{l}^{i j}=\frac{\partial E_{l}}{\partial A_{i j}}$, then we have
    \begin{align}
     \sum_{i, j} \dot{E}_{l}^{i j} g_{ij} &=l E_{l-1} ,\label{2.6}\\
     \sum_{i, j} \dot{E}_{l}^{i j} A_{i j} &=l E_{l} ,\label{2.7}\\
     \sum_{i, j} \dot{E}_{l}^{i j}\left(A^{2}\right)_{i j} &=n E_{1} E_{l}-(n-l) E_{l+1}, \label{2.8}
     \end{align}
    where $\left(A^{2}\right)_{i j}=\sum_{l=1}^{n} A_{i l} A_{l j}$.
\end{lemma}

Before reviewing the famous Minkowski formula and the Newton-Maclaurin inequality in $\mathbb{H}^{n+1}$, we introduce the definition of strictly $k$-convex hypersurfaces.

\begin{definition}\label{def1.14}(\cite{GL20})
    For a bounded domain $\Omega \subset \mathbb{H}^{n+1}$, it is called strictly $k$-convex if the principal curvature $\kappa=(\kappa_{1},\cdots,\kappa_{n})\in \Gamma^{+}_{k}$, where $\Gamma^{+}_{k}$ is the Garding cone 
    \begin{align*}
        \Gamma^{+}_{k}=\{\kappa \in \mathbb{R}^{n} |E_{i}(\kappa)>0,i=1,\cdots,k \}.
    \end{align*}
\end{definition}

\begin{lemma}(\cite{G2013})\label{lem2.4}
    If $\kappa \in \Gamma_{m}^{+}$, the following inequality is called the Newton-MacLaurin inequality
    \begin{align}\label{2.10}
        E_{m+1}(\kappa) E_{l-1}(\kappa) \leq E_{l}(\kappa) E_{m}(\kappa), \quad 1 \leq l \leq m.
    \end{align}
    Equality holds if and only if $\kappa_{1}=\cdots=\kappa_{n}$.
\end{lemma}

\begin{lemma} (\cite{{GL2015}})\label{lem2.3}
    Let $\Sigma$ be a smooth closed hypersurface in $\mathbb{H}^{n+1}$. Then
    \begin{align}\label{2.9}
        \int_{\Sigma} \lambda^{'}E_{l-1}(\kappa) d\mu= \int_{\Sigma} uE_{l}(\kappa) d\mu.
    \end{align}
\end{lemma}

\subsection{ Evolution equations}
\
\vglue-10pt
 \indent

A family of starshaped hypersurfaces $\Sigma_{t}=X(\Sigma,t)$ satisfies the
general flow as below
\begin{align}\label{2.11}
    \begin{cases}
        \frac{\partial}{\partial t}X(x,t)=F(x,t) \nu(x,t),
       \\
       X(\cdot,0)=X_{0},
    \end{cases}
\end{align}
where $F(x,t)$ and $\nu(x,t)$ are the velocity function and the unit outer normal vector of $\Sigma_{t}$, respectively. Then, the radial function $r(\xi,t)$ satisfies the following equation (see e.g.,\cite{HLW20})
\begin{align}\label{2.12}
    \begin{cases}
        \frac{\partial}{\partial t}r = F(x,t)v, \qquad \text{on}\quad\mathbb{S}^{n} \times \mathbb{R}^{+},
       \\
       r(\cdot,0)=r_{0}.
    \end{cases}
\end{align}
Moreover, $\varphi(\xi)=\chi(r(\xi))$ satisfies the initial value problem 
\begin{align}\label{2.13}
    \begin{cases}
        \frac{\partial}{\partial t}\varphi = F(x,t) \frac{v}{\lambda}, \qquad \text{on}\quad\mathbb{S}^{n} \times \mathbb{R}^{+},
       \\
       \varphi(\cdot,0)=\varphi_{0}.
    \end{cases}
\end{align}

We also have the following evolution equations for geometric quantities.
\begin{lemma}(\cite{HLX14})\label{lem2.5}
    \begin{align}
        \frac{\partial}{\partial t}g_{ij}&=2F h_{ij},\label{2.14}\\
        \frac{\partial}{\partial t}d\mu_{t}&=nE_{1}F d\mu_{t},\label{2.15}\\
        \frac{\partial}{\partial t}h_{ij}&=-\nabla_{j}\nabla_{i}F+F\left((h^{2})_{ij}+g_{ij}\right),\label{2.16}\\
        \frac{\partial}{\partial t}h^{j}_{i}&=-\nabla^{j}\nabla_{i}F-F\left((h^{2})^{j}_{i}+\delta^{j}_{i}\right),\label{2.17}\\
        \frac{\partial}{\partial t}E_{l-1}&=\frac{\partial E_{l-1}}{\partial h^{j}_{i}}\frac{\partial h^{j}_{i}}{\partial t}=\dot{E}^{ij}_{l-1}\left(-\nabla_{j}\nabla_{i}F - F(h^{2})_{ij}+Fg_{ij} \right) ,\label{2.28}\\
        \frac{\partial}{\partial t}\lambda^{'}&=\left\langle \bar{\nabla} \lambda^{'},\partial_{t}\right\rangle =uF ,\label{2.19} \\
        \frac{\partial}{\partial t}W^{\lambda^{'}}_{0}(\Omega_{t})&=(n+1)\int_{M_{t}}(n+1)\lambda^{'} F d\mu_{t} ,  \label{2.20}
    \end{align}
    where $\nabla$ denotes the Levi-Civita connection on $(\Sigma_{t},g)$ .
\end{lemma}

\section{ Existence and smooth convergence of the flow (\ref{1.9})}

\noindent  It is known that the locally constrained mean curvature flow (\ref{1.9}) can be represented as the PDE of $\varphi (\xi,t)$, i.e., the following scalar equation (\ref{3.3}). In order to obtain the longtime existence of the flow (\ref{1.9}), we first derive a priori estimate for the equation (\ref{3.3}) and then use the standard theory of parabolic partial differential equations. The convergence of the flow (\ref{1.9}) is obtained by a detailed estimation of the gradient of $\varphi$. 

\subsection{ Longtime exists of the flow (\ref{1.9})}
\
\vglue-10pt
 \indent

From (\ref{2.12}), the flow (\ref{1.9}) can be parameterized into a scalar PDE of $r$ as following
\begin{align}\label{3.1}
    \begin{cases}
        \frac{\partial}{\partial t} r =-\Phi_{1}H-\frac{n}{n-1} \frac{\partial \Phi_{1}}{\partial r}\sqrt{1+\lambda^{-2}|D r|^{2}}, \qquad \text{on}\quad\mathbb{S}^{n} \times \mathbb{R}^{+},
        \\
        r(\cdot,0)=r_{0}.
    \end{cases}
\end{align}
There are different representations of $f$ on starshaped hypersurfaces $\Sigma_{t}$, such as
\begin{align*}
    &f:=\Phi_{1}(r(\xi)):\Sigma_{t} \to \mathbb{S}^{n} \to \mathbb{R}; \quad
    r: \mathbb{S}^{n} \to \mathbb{R},\\
    &f:=\Phi_{1}(\varphi(\xi)):\Sigma_{t} \to \mathbb{S}^{n} \to \mathbb{R}; \quad
    \varphi :\mathbb{S}^{n} \to \mathbb{R},
\end{align*}
where $\varphi(\xi) = \chi(r(\xi))$ and $\frac{\partial }{\partial r}\chi=\frac{1}{\lambda(r)}$. Thus, we have 
\begin{align}\label{3.2}
    \frac{\partial \Phi_{1}}{\partial r}=\frac{\partial \Phi_{1}}{\partial \varphi} \frac{1}{\lambda}.
\end{align}
Furthermore, by (\ref{2.13}) and (\ref{3.2}), $\varphi$ satisfies the initial value problem 
\begin{align}\label{3.3}
    \begin{cases}
        \frac{\partial}{\partial t}\varphi=-\Phi_{1}\frac{H}{\lambda}-\frac{n}{n-1} \frac{\partial \Phi_{1}}{\partial \varphi}\frac{\sqrt{1+|D \varphi|^{2}}}{\lambda^{2}}, \qquad \text{on}\quad\mathbb{S}^{n} \times \mathbb{R}^{+},
        \\
        \varphi(\cdot,0)=\varphi_{0}.
    \end{cases}
\end{align}

Firstly, we perform $C^{0}$ estimate of $\varphi$.
\begin{proposition}\label{prop3.1}
    Let $\varphi \in C^{\infty}\left(\mathbb{S}^{n} \times [0,T) \right)$ be a solution to the initial value problem (\ref{3.3}), then there are positive constants $C_{1}=C_{1}( \varepsilon, \varphi_{min}(0))$ and $C_{2}=C_{2}( \delta, \varphi_{max}(0))$, such that
    \begin{align}\label{3.4}
          C_{1} \leq \varphi(\cdot,t) \leq C_{2}.
    \end{align}
\end{proposition}
\noindent{\it \bf Proof.}~~ Let $\varphi_{min}(t):=\underset{\xi \in
\mathbb{S}^{n}}{\min}  \varphi(\cdot,t)$, then $D \varphi_{min} =0$ and $D^{2} \varphi_{min} \ge 0$. In view of (\ref{2.5}),  at the spatial minimum point of $\varphi$, we have 
\begin{align*}
    H(\varphi_{min})&=n\frac{\lambda^{'}}{\lambda}-\frac{1}{\lambda}D^{2}\varphi_{min}.
\end{align*}
From (\ref{3.3}) and the above equality, we get 
\begin{align*}
    \frac{\partial }{\partial t}\varphi_{min}&=-\lambda^{-2}\left(n\Phi_{1}\lambda^{'}+\frac{n}{n-1}\frac{\partial \Phi_{1}}{\partial \varphi}\right)+\lambda^{-2} \Phi_{1} D^{2}\varphi_{min}\\
    &\geq -\lambda^{-2}\left(n\Phi_{1}\lambda^{'}+\frac{n}{n-1}\frac{\partial \Phi_{1}}{\partial \varphi}\right).
\end{align*}
Using (\ref{3.2}), $\widehat{\Phi}_{1}(r)=n\Phi_{1}\frac{\lambda^{'}}{\lambda^{2}}+\frac{n}{n-1} \frac{\partial \Phi_{1}}{\partial r}\frac{1}{\lambda}$ can be expressed as $\widehat{\Phi}_{1}(\varphi)=n\Phi_{1}\frac{\lambda^{'}}{\lambda^{2}}+\frac{n}{n-1}\frac{\partial \Phi_{1}}{\partial \varphi}\frac{1}{\lambda^{2}}$. Thus
\begin{align*}
    \frac{\partial }{\partial t}\varphi_{min} \geq - \widehat{\Phi}_{1}(\varphi_{min}).
\end{align*}
If $\widehat{\Phi}_{1}(\varphi_{min}) \leq 0$, we have $\frac{\partial }{\partial t}\varphi_{min} \geq 0$. Then
\begin{align*}
    \varphi_{min}(t) \geq \varphi_{min}(0).
\end{align*}
If $\widehat{\Phi}_{1}(\varphi_{min}) > 0$, then, according to Assumption \ref{as1.10}, there exists a zero point for $\widehat{\Phi}_{1}(\varphi)$. Consequently, there exist constants $\varepsilon$ and $ \delta$ such that $0<\varepsilon<\delta$ and 
\begin{align*}
    \widehat{\Phi}_{1}(\varepsilon)<0,\qquad \widehat{\Phi}_{1}(\delta)>0.
\end{align*}
Thus $\widehat{\Phi}_{1}(\varphi_{min}) > \widehat{\Phi}_{1}(\varepsilon)$. Moreover, Since $\widehat{\Phi}_{1}(r)$ is monotonically increasing with respect to $r$, it follows from (\ref{3.2}) that $\widehat{\Phi}_{1}(\varphi)$ is monotonically increasing with respect to $\varphi$. Then, we can deduce that
\begin{align*}
    \varphi_{min}(t) > \varepsilon.
\end{align*}
Combining the above two cases, we get
\begin{align*}
    \varphi_{min}(t) \geq C_{1}:= \min \{  \varepsilon,\varphi_{min}(0) \}.
\end{align*}

In the same way, let $\varphi_{max}(t):=\underset{\xi \in \mathbb{S}^{n}}{\max}  \varphi(\cdot,t)$, we have $D \varphi_{max} =0$, $D^{2} \varphi_{max} \leq 0$ and 
\begin{align*}
    \frac{\partial }{\partial t} \varphi_{max} \leq -\lambda^{-2}\left(n\Phi_{1}\lambda^{'}+\frac{n}{n-1}\frac{\partial \Phi_{1}}{\partial \varphi}\right) = -\widehat{\Phi}_{1}(\varphi_{max}).
\end{align*}
If $\widehat{\Phi}_{1}(\varphi_{max})\geq 0$, we obtain $\frac{\partial }{\partial t} \varphi_{max} \leq 0$, i.e., $\varphi_{max}(t)\leq \varphi_{max}(0)$. If $\widehat{\Phi}_{1}(\varphi_{max})< 0$, then $\widehat{\Phi}_{1}(\varphi_{max}) < \widehat{\Phi}_{1}(\delta)$, i.e., $\varphi_{max}(t) < \delta$. Combining the above two situations, we obtain
\begin{align*}
    \varphi_{max}(t) \leq C_{2}:= \min \{\delta, \varphi_{max}(0)\}.
\end{align*}
This completes the proof of Proposition \ref{prop3.1}.\hfill${\square}$

In the following, we prove the $C^{1}$ estimate of $\varphi$.
\begin{proposition}\label{prop3.2}
    Let $\varphi \in {C^{\infty}(\mathbb{S}^{n} \times [0,T))}$ be a solution to the initial value problem (\ref{3.3}). For any time $t \in [0,T)$, there is a positive constant $C$ depending on $\Sigma_{0}$, such that
    \begin{align}\label{3.5}
        \underset{\xi \in \mathbb{S}^{n}}{\max} |D \varphi(\cdot,t)|  \leq C.
    \end{align}
\end{proposition}
\noindent{\it \bf Proof.}~~Let $\psi =\frac{1}{2}|D \varphi|^{2}$, we have $\frac{\partial}{\partial t}\psi=D^{k}\varphi \cdot D_{k}\left(\frac{\partial }{\partial t}\varphi\right)=D^{k}\varphi\cdot \frac{\partial}{\partial t}\left(D_{k}\varphi\right)$. According to (\ref{3.3}), the evolution equation of $\psi$ is 
\begin{align}\label{3.6}
    \frac{\partial }{\partial t}\psi=
    -D^{k}\varphi \left[D_{k}(\Phi_{1})\frac{H}{\lambda}+\Phi_{1} D_{k}\left(\frac{H}{\lambda}\right)+\frac{n}{n-1}D_{k}\left(\frac{\partial \Phi_{1}}{\partial \varphi}\right)\frac{v}{\lambda^{2}}+\frac{n}{n-1}\frac{\partial \Phi_{1}}{\partial \varphi}D_{k}\left(\frac{v}{\lambda^{2}}\right)\right],
\end{align}
where
\begin{align*}
    D_{k}\Phi_{1}=\frac{\partial \Phi_{1}}{\partial \varphi}D_{k}\varphi, \qquad D_{k}\left(\frac{\partial \Phi_{1}}{\partial \varphi}\right)=\frac{\partial^{2} \Phi_{1}}{\partial \varphi^{2}}D_{k}\varphi.
\end{align*}
Suppose $\psi$ attains the spatial maximum at point $\left(\xi_{t},t\right) \in \left(\mathbb{S}^{n} \times [0,T)\right)$, $\psi_{max} := \psi (\xi_{t},t)=\underset{\xi \in \mathbb{S}^{n}}{\max}\psi(\xi,t) $, then 
\begin{align}
    D \psi_{max}&=D^{m}\varphi D_{km}\varphi =0,    \label{3.7}\\
    D^{2}\psi_{max}&=D^{km}\varphi D_{km}\varphi + D^{k}\varphi D_{klm}\varphi  \leq 0,  \qquad k=1,\dots ,n, \label{3.8}
\end{align}
and at the spatial maximum point $\left(\xi_{t},t\right)$ of $\psi$, there also holds
\begin{align*}
    D_{k}v=D_{k}\left(\sqrt{1+|D \varphi|^{2}}\right)=&-v^{-1}D^{m}\varphi D_{km}\varphi =0 ,\\
    D^{k}\varphi \cdot H  =& n (v\lambda)^{-1}\lambda^{'} D^{k}\varphi,\\
    D_{k}H=& n(v\lambda)^{-1} D_{k}\varphi  \left(\lambda^{2} -(\lambda^{'})^{2}\right) -(v \lambda)^{-1}\left(\sigma^{ki}-\frac{\varphi^{k}\varphi^{i}}{v^{2}}\right) D_{k}\varphi_{ki},\\
    D^{k}\varphi \cdot {D_{k}H} \geq & n(v\lambda)^{-1}|D \varphi|^{2} \left(\lambda^{2}-(\lambda^{'})^{2}\right).
\end{align*}
Substituting  (\ref{3.7}), (\ref{3.8}) and the above equations into (\ref{3.6}), we obtain
\begin{align}\label{3.9}
    \frac{\partial }{\partial t}\psi_{max} \leq & 2v^{-1}\psi_{max} \left[\frac{2n}{n-1}\frac{\partial \Phi_{1}}{\partial \varphi}\frac{\lambda^{'}}{\lambda^{2}}v^{2}-\frac{n}{n-1}\frac{\partial^{2}\Phi_{1}}{\partial \varphi^{2}}\frac{1}{\lambda^{2}}v^{2}-n\frac{\partial \Phi_{1}}{\partial \varphi}\frac{\lambda^{'}}{\lambda^{2}}-n\Phi_{1}+2n\Phi_{1} \frac{(\lambda^{'})^{2}}{\lambda^{2}}\right]  \notag \\
    =& 2v^{-1}\psi_{max} \left[\frac{2n}{n-1}\frac{\partial \Phi_{1}}{\partial \varphi}\frac{\lambda^{'}}{\lambda^{2}}-\frac{n}{n-1}\frac{\partial^{2}\Phi_{1}}{\partial \varphi^{2}}\frac{1}{\lambda^{2}}-n\frac{\partial \Phi_{1}}{\partial \varphi}\frac{\lambda^{'}}{\lambda^{2}}-n\Phi_{1}+2n\Phi_{1}\frac{(\lambda^{'})^{2}}{\lambda^{2}}\right]\\
    &+ 4v^{-1}\psi^{2}_{max} \left[\frac{2n}{n-1}\frac{\partial \Phi_{1}}{\partial \varphi}\frac{\lambda^{'}}{\lambda^{2}}-\frac{n}{n-1}\frac{\partial^{2} \Phi_{1}}{\partial \varphi^{2}}\frac{1}{\lambda^{2}}\right].   \notag  
\end{align}
Denote by
\begin{align*}
    C_{3}(\varphi(t))&=2v^{-1}\left[\frac{2n}{n-1}\frac{\partial \Phi_{1}}{\partial \varphi}\frac{\lambda^{'}}{\lambda^{2}}-\frac{n}{n-1}\frac{\partial^{2}\Phi_{1}}{\partial \varphi^{2}}\frac{1}{\lambda^{2}}-n\frac{\partial \Phi_{1}}{\partial \varphi}\frac{\lambda^{'}}{\lambda^{2}}-n\Phi_{1}+2n\Phi_{1}\frac{(\lambda^{'})^{2}}{\lambda^{2}}\right], \\
    C_{4}(\varphi(t))&=4v^{-1} \left[\frac{2n}{n-1}\frac{\partial \Phi_{1}}{\partial \varphi}\frac{\lambda^{'}}{\lambda^{2}}-\frac{n}{n-1}\frac{\partial^{2} \Phi_{1}}{\partial \varphi^{2}}\frac{1}{\lambda^{2}}\right] .
\end{align*}
Thus, (\ref{3.9}) can be expressed as 
\begin{align*}
    \frac{\partial }{\partial t}\psi_{max} \leq C_{3}(\varphi(t))\psi_{max} + C_{4}(\varphi(t))\psi^{2}_{max}.
\end{align*}

We now prove that $C_{3}(\varphi(t)) \leq 0$ and $C_{4}(\varphi(t)) < 0$. From Assumption \ref{as1.10}, $\widehat{\Phi}_{1}(r)$ is monotonically increasing with respect to $r$ and 
\begin{align*}
    \frac{\partial \widehat{\Phi}_{1}}{ \partial r}  = \left(\frac{n}{\lambda} - \frac{2n(\lambda^{'})^{2}}{\lambda^{3}}\right)\Phi_{1} +\left(\frac{n\lambda^{'}}{\lambda^{2}} - \frac{n}{n-1}\frac{\lambda^{'}}{\lambda^{2}}\right)\frac{\partial \Phi_{1}}{\partial r} + \frac{n}{n-1}\frac{1}{\lambda}\frac{\partial^{2} \Phi_{1}}{\partial r^{2}}. 
\end{align*}
Then
\begin{align}\label{3.10}
    \frac{n}{n-1}\frac{\partial \Phi_{1}}{\partial r}\frac{\lambda^{'}}{\lambda^{2}}-\frac{n}{n-1}\frac{\partial^{2} \Phi_{1}}{\partial r^{2}}\frac{1}{\lambda}\leq n\frac{\partial \Phi_{1}}{\partial r}\frac{\lambda^{'}}{\lambda^{2}}+n\Phi_{1}\frac{1}{\lambda}-2n\Phi_{1}\frac{(\lambda^{'})^{2}}{\lambda^{3}}.
\end{align}
Combining $\frac{\partial \Phi_{1}}{\partial r} = \frac{\partial \Phi_{1}}{\partial \varphi}\frac{1}{\lambda}$ and $\frac{\partial \varphi}{\partial r} = \frac{1}{\lambda}$, we get 
\begin{align*}
    \frac{\partial^{2} \Phi_{1}}{\partial r^{2}} = \frac{\partial^{2} \Phi_{1}}{\partial \varphi^{2}}\frac{1}{\lambda^{2}} - \frac{\partial \Phi_{1}}{\partial \varphi}\frac{\lambda^{'}}{\lambda^{2}}.
\end{align*}
Thus, the inequality (\ref{3.10}) is equivalent to
\begin{align}\label{3.11}
    \frac{2n}{n-1}\frac{\partial \Phi_{1}}{\partial \varphi}\frac{\lambda^{'}}{\lambda^{2}}-\frac{n}{n-1}\frac{\partial^{2} \Phi_{1}}{\partial \varphi^{2}}\frac{1}{\lambda^{2}}\leq n\frac{\partial \Phi_{1}}{\partial \varphi}\frac{\lambda^{'}}{\lambda^{2}}+n\Phi_{1}-2n\Phi_{1}\frac{(\lambda^{'})^{2}}{\lambda^{2}}.
\end{align}
By (\ref{3.11}), it is not difficult to find that $C_{3}(\varphi(t)) \leq 0$.

Next, we prove $C_{4}(\varphi(t)) < 0$. Suppose that 
\begin{align*}
    C_{4}(\varphi(t)) =4v^{-1} \left[\frac{2n}{n-1}\frac{\partial \Phi_{1}}{\partial \varphi}\frac{\lambda^{'}}{\lambda^{2}}-\frac{n}{n-1}\frac{\partial^{2} \Phi_{1}}{\partial \varphi^{2}}\frac{1}{\lambda^{2}}\right] \geq 0,
\end{align*}
then $\frac{2n}{n-1}\frac{\partial \Phi_{1}}{\partial \varphi}\frac{\lambda^{'}}{\lambda^{2}}-\frac{n}{n-1}\frac{\partial^{2} \Phi_{1}}{\partial \varphi^{2}}\frac{1}{\lambda^{2}} \geq  0$. From Proposition \ref{prop3.1}, we have $C_{1} \leq \varphi(\xi ,t) \leq C_{2}$ holds for all $(\xi,t) \in \mathbb{S}^{n} \times [0,T)$, where  $C_{1}= min\{\varepsilon, \varphi_{min}(0)\}$ and $C_{1},C_{2}, \varepsilon \in \mathbb{R}^{+}$.  Integrating both sides of $\frac{2n}{n-1}\frac{\partial \Phi_{1}}{\partial \varphi}\frac{\lambda^{'}}{\lambda^{2}}-\frac{n}{n-1}\frac{\partial^{2} \Phi_{1}}{\partial \varphi^{2}}\frac{1}{\lambda^{2}} \geq  0$ yields
\begin{align*}
    0 &\leq  \frac{n}{n-1} \int_{ C_{1}}^{\varphi} \left(2\frac{\partial \Phi_{1}}{\partial z}\frac{\lambda^{'}}{\lambda^{2}} - \frac{\partial^{2}\Phi_{1}}{\partial z^{2}}\frac{1}{\lambda^{2}} \right) dz = -\frac{n}{n-1}\int_{C_{1}}^{\varphi} \frac{\partial }{\partial z} \left(\frac{1}{\lambda^{2}}\frac{\partial \Phi_{1}}{\partial z}\right) d z
    \\
    &= \frac{n}{n-1} \left[\left(\frac{1}{\lambda^{2}}\frac{\partial \Phi_{1}}{\partial \varphi}(C_{1})\right) - \left(\frac{1}{\lambda^{2}}\frac{\partial \Phi_{1}}{\partial \varphi}\right) \right],
\end{align*}
where we used $\frac{\partial \lambda}{\partial \varphi} = \lambda^{'}\lambda$. Thus,
\begin{align*}
    \frac{1}{\lambda^{2}}\frac{\partial \Phi_{1}}{\partial \varphi} \leq \frac{1}{\lambda^{2}}\frac{\partial \Phi_{1}}{\partial \varphi}(C_{1}).
\end{align*}

In fact, $\frac{1}{\lambda^{2}}\frac{\partial \Phi_{1}}{\partial \varphi}(C_{1}) < 0$.  To prove this assertion, we need to consider two cases: (1) when $C_{1} = \varepsilon $, it follows from the proof of Proposition \ref{prop3.1} that $\widehat{\Phi}_{1}(\varepsilon) < 0$;  (2) when $C_{1} =\varphi_{min}(0) $, we have  $\varphi_{min}(0) \leq \varepsilon$.  From (\ref{3.11}), we know that $\widehat{\Phi}_{1}(\varphi)$ is monotonically increasing with respect to $\varphi$, then $\widehat{\Phi}_{1}(\varphi_{min}(0)) \leq \widehat{\Phi}_{1}(\varepsilon) < 0$. Thus, $\widehat{\Phi}_{1}(C_{1}) <0$.  Since $\widehat{\Phi}_{1}(\varphi)= \frac{n}{n-1} \frac{\partial \Phi_{1}}{\partial \varphi}\frac{1}{\lambda^{2}}+n\Phi_{1}\frac{\lambda^{'}}{\lambda^{2}}$, it follows that
\begin{align*}
    \widehat{\Phi}_{1}(C_{1}) =  \frac{n}{n-1} \frac{\partial \Phi_{1}}{\partial \varphi}(C_{1})\frac{1}{\lambda^{2}}+n\Phi_{1}(C_{1})\frac{\lambda^{'}}{\lambda^{2}}< 0,
\end{align*}
and 
\begin{align*}
    \frac{1}{\lambda^{2}}\frac{\partial \Phi_{1}}{\partial \varphi}(C_{1}) \leq -(n-1)\frac{\lambda^{'}}{\lambda^{2}}\Phi_{1}(C_{1})< 0 .
\end{align*}
Thus, 
\begin{align*}
    \frac{1}{\lambda^{2}}\frac{\partial \Phi_{1}}{\partial \varphi}  \leq \frac{1}{\lambda^{2}}\frac{\partial \Phi_{1}}{\partial \varphi} (C_{1})< 0
\end{align*}
holds for any $\varphi \in [C_{1},C_{2}]$, and 
\begin{align*}
    n\frac{\partial \Phi_{1}}{\partial \varphi}\frac{\lambda^{'}}{\lambda^{2}}+n\Phi_{1}-2n\Phi_{1}\frac{(\lambda^{'})^{2}}{\lambda^{2}} <0.
\end{align*}
Combining (\ref{3.11}) and the above inequality, we get 
\begin{align*}
    \frac{2n}{n-1}\frac{\partial \Phi_{1}}{\partial \varphi}\frac{\lambda^{'}}{\lambda^{2}}-\frac{n}{n-1}\frac{\partial^{2} \Phi_{1}}{\partial \varphi^{2}}\frac{1}{\lambda^{2}}\leq n\frac{\partial \Phi_{1}}{\partial \varphi}\frac{\lambda^{'}}{\lambda^{2}}+n\Phi_{1}-2n\Phi_{1}\frac{(\lambda^{'})^{2}}{\lambda^{2}} < 0,
\end{align*}
which contradicts $\frac{2n}{n-1}\frac{\partial \Phi_{1}}{\partial \varphi}\frac{\lambda^{'}}{\lambda^{2}}-\frac{n}{n-1}\frac{\partial^{2} \Phi_{1}}{\partial \varphi^{2}}\frac{1}{\lambda^{2}} \geq 0$. Thus, we have $\frac{2n}{n-1}\frac{\partial \Phi_{1}}{\partial \varphi}\frac{\lambda^{'}}{\lambda^{2}}-\frac{n}{n-1}\frac{\partial^{2} \Phi_{1}}{\partial \varphi^{2}}\frac{1}{\lambda^{2}} < 0$, i.e., $C_{4}(\varphi(t)) < 0$.

From $C_{3}(\varphi(t)) \leq 0 $ and  $C_{4}(\varphi(t)) < 0 $, the inequality (\ref{3.9}) can be scaled as $\frac{\partial }{\partial t} \psi_{max}\leq  0$. Hence, $\psi_{max} (t) \leq \psi_{max}(0) $, which implies (\ref{3.5}). This completes the proof of Proposition \ref{prop3.2}.

\hfill${\square}$

Equation (\ref{3.3}) can be rewritten as a scalar parabolic PDE in divergent form as follows
\begin{align*}
    \frac{\partial \varphi}{\partial t}&=\frac{\Phi_{1}}{\lambda^{2}v}\left(\sigma^{ki}-\frac{\varphi^{k}\varphi^{i}}{v^{2}}\right)\varphi_{ik}-n\lambda^{'}\frac{\Phi_{1}}{\lambda^{2}v}-\frac{n}{n-1}\frac{\partial \Phi_{1}}{\partial \varphi}\frac{v}{\lambda^{2}} \\
    &=\operatorname{div}\left(\frac{\Phi_{1}}{\lambda^{2}v}\cdot D \varphi\right)+\frac{|D \varphi|^{2}}{\lambda^{2}v}\left(2\Phi_{1}\lambda^{'}-\frac{\partial \Phi_{1}}{\partial \varphi}\right)-n\Phi_{1}\frac{\lambda^{'}}{\lambda^{2}v}-\frac{n}{n-1}\frac{\partial \Phi_{1}}{\partial \varphi}\frac{v}{\lambda^{2}},
\end{align*}
where we used $v_{i}=\frac{\varphi^{j}\varphi_{ji}}{v}$, $D_{i}\lambda=\lambda^{'}D_{i}r=\lambda \lambda^{'} D_{i}\varphi$, and $D_{i}\lambda^{'}=\lambda D_{i}r=\lambda^{2}D_{i}\varphi$. According to the classical theory of parabolic PDEs in divergent form \cite{LSU}, the higher regularity estimate for the solution $\varphi$ can be derived from the uniform $C^{1}$ estimate in Proposition \ref{prop3.2}. Hence, the solution $\Sigma_{t}$ of the flow (\ref{1.9}) has longtime existence.

\subsection{ Exponential convergence of the flow (\ref{1.9})}
\
\vglue-10pt
 \indent

Finally, we prove the flow (\ref{1.9}) exponentially converges to a geodesic sphere by making an exact estimate of $|D\varphi|$.
\begin{proposition}\label{prop3.3}
    Let $\varphi \in {C^{\infty}(\mathbb{S}^{n} \times [0,\infty))}$ be a solution to the initial value problem (\ref{3.3}). For any time $t \in [0,\infty)$, there exist two positive constants $\overline{C}$ and $\gamma$ that depend only on $\Sigma_{0}$, such that
    \begin{align}\label{3.12}
        \underset{\xi\in \mathbb{S}^{n}}{\max} |D \varphi(\cdot, t)|^{2}\leq \overline{C}e^{-\gamma t}.
    \end{align}
\end{proposition}
\noindent{\it \bf Proof.}~~From the proof of Proposition \ref{prop3.2}, we have, at the spatial maximum point of $\psi$,
\begin{align*}
    \frac{\partial }{\partial t}\psi_{max} \leq C_{3}(\varphi(t))\psi_{max} + C_{4}(\varphi(t))\psi^{2}_{max},
\end{align*}
and $C_{3}(\varphi(t)) \leq 0$, $C_{4}(\varphi(t)) < 0$. Since $\Phi_{1}(r)$ is a smooth function, we have $\frac{\partial \Phi_{1}}{\partial r} = \lambda^{-1} \frac{\partial \Phi_{1}}{\partial \varphi}$ and $\frac{\partial^{2} \Phi_{1}}{\partial r^{2}}=\lambda^{-2}\left(\frac{\partial^{2} \Phi_{1}}{\partial \varphi^{2}}-\lambda^{'}\frac{\partial \Phi_{1}}{\partial \varphi} \right)$ are continuous functions. Thus, given that $C_{1} \leq \varphi(t) \leq C_{2}$, it follows that $ C_{3}(\varphi(t)) $ and $ C_{4}(\varphi(t)) $ are uniformly bounded continuous functions.

It is easy to show that there exist positive constants $\gamma$ and $C$ such that $\psi (t) \leq Ce^{-2\gamma t}$. In fact, for $C_{3}(\varphi(t)) < 0$, we denote $-2\gamma:=  \max_{\varphi (t) \in [C_{1},C_{2}]} C_{3}(\varphi (t)) $, then $\gamma$ is a positive constant and
\begin{align*}
    \frac{\partial }{\partial t}\psi_{max}  \leq C_{3}(\varphi(t)) \psi_{max} \leq -2\gamma \psi_{max}.
\end{align*}
Integrating both ends of the above inequality yields
\begin{align*}
    \psi_{max}(t) \leq Ce^{-2\gamma t},
\end{align*}
where $C$ is a positive constant that depends only on $\Sigma_{0}$.

For $C_{3}(\varphi(t))=2v^{-1}\left[\frac{2n}{n-1}\frac{\partial \Phi_{1}}{\partial \varphi}\frac{\lambda^{'}}{\lambda^{2}}-\frac{n}{n-1}\frac{\partial^{2}\Phi_{1}}{\partial \varphi^{2}}\frac{1}{\lambda^{2}}-n\frac{\partial \Phi_{1}}{\partial \varphi}\frac{\lambda^{'}}{\lambda^{2}}-n\Phi_{1}+2n\Phi_{1}\frac{(\lambda^{'})^{2}}{\lambda^{2}}\right] = 0$, we have 
\begin{align*}
    \frac{3n-n^{2}}{n-1}\frac{\partial \Phi_{1}}{\partial \varphi}\frac{\lambda^{'}}{\lambda^{2}}-\frac{n}{n-1}\frac{\partial^{2} \Phi_{1}}{\partial \varphi^{2}}\frac{1}{\lambda^{2}}+2n\Phi_{1}\frac{(\lambda^{'})^{2}}{\lambda^{2}}-n\Phi_{1}=0,
\end{align*}
which is equivalent to 
\begin{align*}
    \frac{2n-n^{2}}{n-1}\frac{\partial \Phi_{1}}{\partial r} \frac{\lambda^{'}}{\lambda}-\frac{n}{n-1}\frac{\partial^{2} \Phi_{1}}{\partial r^{2}} + 2n\Phi_{1}\left(\frac{\lambda^{'}}{\lambda}\right)^{2} -n\Phi_{1} = 0.
\end{align*}
Solving the above ODE, we can obtain the general solution as
\begin{align*}
    \Phi_{1}(r) =C_{5}\lambda^{1-n}(r) = C_{5}({\rm{sinh}}r)^{1-n},
\end{align*}
where $C_{5}$ is a positive constant. Since $\widehat{\Phi}_{1}(r)=n\Phi_{1}\frac{\lambda^{'}}{\lambda^{2}}+ \frac{n}{n-1} \frac{\partial \Phi_{1}}{\partial r}\frac{1}{\lambda}$, it follows that 
\begin{align*}
    0 \equiv \widehat{\Phi}_{1}(r)  = \widehat{\Phi}_{1}(\varphi),  \qquad \text{ for any }  \varphi(t) \in [C_{1},C_{2}],
\end{align*}
which does not satisfy Assumption \ref{as1.10}. This completes the proof of Proposition \ref{prop3.3}.
\hfill${\square}$

\noindent {\bf Proof of Theorem \ref{th1.11}}: It can be inferred from Proposition \ref{prop3.3} that there exists a sequence of times $\{ t_{i} \}$ such that $\lim_{t_{i} \to \infty}
r_{t_{i}}=r_{\infty}$.  By  the interpolation inequality and Sobolev embedding theorem on $\mathbb{S}^{n}$, one gets the convergence of $r_{t_{i}}$ to $r_{\infty}$ in $C^\infty$-topology. Since the velocity of the flow (\ref{1.9}) dose not contain global terms, by the comparison principle, the radius $r_{\infty}$ is unique.
Moreover, $B_{r_{\infty}}$ must be centered at the origin. Otherwise, there will be a conflict with the proof of Proposition \ref{prop3.1}. 
\hfill${\square}$

\section{ Michael-Simon type inequalities for $k$-th mean curvatures}

\noindent In this section, we will use the smooth convergence results of the flows (\ref{1.9}) and (\ref{1.12}) to give the proofs of Theorem \ref{th1.12} and Theorem \ref{th1.15} respectively.

\subsection{ Sharp Michael-Simon type inequality for
mean curvature}
\
\vglue-10pt
 \indent

In this section we will prove the desired geometric inequalities (\ref{1.10}) and (\ref{1.11}) for starshaped hypersurfaces.

\noindent {\bf Proof of Theorem \ref{th1.12}}: First, we reduce the inequality (\ref{1.10}) by scaling. Assume that
\begin{align*}
    \int_{M_{0}}  \lambda^{'} \sqrt{f^{2} E_{1}^{2}+|\nabla^{M_{0}} f|^{2}}-\int_{M_{0}}\left\langle \bar{\nabla}\left(f\lambda^{'}\right),\nu \right\rangle  +\int_{\partial M_{0}}  f
    = \int_{M_{0}}f^{\frac{n}{n-1}}.
\end{align*}
This normalization guarantees that we can find the solution
$\vartheta: M_{0} \rightarrow \mathbb{R}$ to the following equation
\begin{align*}
    \operatorname{div}_{M_{0}}\left(f \nabla^{M} \vartheta \right)= f^{\frac{n}{n-1}}-\lambda^{'} \sqrt{f^{2} E_{1}^{2}+|\nabla^{M} f|^{2}} + \left\langle \bar{\nabla}\left(f\lambda^{'}\right),\nu \right\rangle
\end{align*}
on $M_{0}$, and $\left\langle\nabla^{M_{0}} \vartheta , \vec{n}\right\rangle=1$ on $\partial M_{0}$. Here, $\vec{n}$ denotes the co-normal to $M_{0}$. According to the standard elliptic regularity theory, $\vartheta$ belongs to $C^{2, \beta}$ for each $0<\beta<1$. Thus, we only need to prove $\int_{M_{0}}f^{\frac{n}{n-1}} d\mu \geq \omega_{n}$. Note that $M_{t} = \Sigma_{t}$ or $ M_{t} = \Sigma_{t} \cup \partial \Sigma_{t}$, $t\in[0,\infty)$, then 
\begin{align}\label{4.1}
    \int_{M_{t}}f^{\frac{n}{n-1}} d\mu \geq \int_{\Sigma_{t}}f^{\frac{n}{n-1}} d\mu \geq \omega_{n},
\end{align}
and it is only necessary to prove the second inequality. 

Secondly, we prove the monotonicity of $\int_{\Sigma_{t}}f^{\frac{n}{n-1}}d\mu_{t}$, which is the key point. Since $f\big|_{\Sigma_{t}} = \Phi_{1} \circ r(\xi,t)$, it follows that $\int_{\Sigma_{t}}f^{\frac{n}{n-1}}d\mu_{t}=\int_{\Sigma_{t}}\Phi_{1}^{\frac{n}{n-1}}(r)d\mu_{t}$.  Along the flow (\ref{1.9}), we have
\begin{align*}
    \frac{\partial}{\partial t}\int_{\Sigma_{t}}\Phi_{1}^{\frac{n}{n-1}}d\mu_{t} &=\int_{\Sigma_{t}}\left(\frac{n}{n-1}\Phi_{1}^{\frac{1}{n-1}}\frac{\partial \Phi_{1}}{\partial t} + \Phi_{1}^{\frac{n}{n-1}}F H \right)d\mu_{t}\\
    &=\int_{\Sigma_{t}}\Phi_{1}^{\frac{1}{n-1}}\left(\frac{n}{n-1}\frac{\partial \Phi_{1}}{\partial r}\sqrt{1+\lambda^{-2}|D r|^{2}}+\Phi_{1} H\right)F d\mu_{t},\\
\end{align*}
where $F=-\left(\frac{n}{n-1}\frac{\partial \Phi_{1}}{\partial r} +\Phi_{1} \frac{H}{\sqrt{1+\lambda^{-2}|D r|^{2}}} \right)$. Thus
\begin{align}\label{4.2}
    \frac{\partial}{\partial t}&\int_{\Sigma_{t}}\Phi_{1}^{\frac{n}{n-1}}d\mu_{t} =-\int_{\Sigma_{t}}\frac{\Phi_{1}^{\frac{1}{n-1}}}{\sqrt{1+\lambda^{-2}|D r|^{2}}}\left(\frac{n}{n-1}\frac{\partial \Phi_{1}}{\partial r}\sqrt{1+\lambda^{-2}|D r|^{2}} +\Phi_{1} H\right)^{2} d\mu_{t} \leq 0.
\end{align}
The monotonically decreasing property of $\int_{\Sigma_{t}}\Phi_{1}^{\frac{n}{n-1}}d\mu_{t}$ yields
\begin{align*}
    \int_{\Sigma_{0}}\Phi_{1}^{\frac{n}{n-1}}d\mu \geq \int_{\Sigma_{t}}\Phi_{1}^{\frac{n}{n-1}}d\mu_{t} \geq \int_{\Sigma_{\infty}}\Phi_{1}^{\frac{n}{n-1}}d\mu_{\infty}.
\end{align*}
From the convergence result of the flow (\ref{1.9}) that $\Sigma_{\infty}=B_{r_{\infty}}$, we have
\begin{align*}
    \nabla r_{\infty} =0, \qquad H |_{\Sigma_{\infty}} = \frac{n\lambda^{'}(r_{\infty})}{\lambda(r_{\infty})},
\end{align*}
and
\begin{align*}
    \nabla \Phi_{1}(r_{\infty}) =\frac{\partial \Phi_{1}}{\partial r_{\infty}} \nabla r_{\infty}= 0.
\end{align*}
Therefore, $f\big|_{\Sigma_{\infty}}=\Phi_{1}(r_{\infty})$ is constant and
\begin{align}\label{4.3}
    \int_{\Sigma_{\infty}}f^{\frac{n}{n-1}}d\mu_{\infty} =\int_{B_{r_{\infty}}}\Phi_{1}(r_{\infty})^{\frac{n}{n-1}}d\mu_{\mathbb{H}^{n+1}}= \Phi^{\frac{n}{n-1}}_{1}(r_{\infty})\lambda^{n}(r_{\infty})\omega_{n}.
\end{align}
The same, $\frac{\partial }{\partial t}\int_{\Sigma_{t}}\Phi_{1}^{\frac{n}{n-1}}d\mu_{t} = 0$ when $t \to \infty$ and by (\ref{4.2}), we have
\begin{align}\label{4.4}
    \frac{n}{n-1}\frac{\partial \Phi_{1}}{\partial r}\sqrt{1+\lambda^{-2}|D r|^{2}}+\Phi_{1} H = 0,
\end{align}
i.e.,
\begin{align*}
    n\Phi_{1}(r_{\infty})\frac{\lambda^{'}(r_{\infty})}{\lambda(r_{\infty})}+\frac{n}{n-1}\frac{\partial \Phi_{1}}{\partial r}(r_{\infty}) = 0.
\end{align*}
It is not difficult to verify that one of the solutions to the above equation is
\begin{align}\label{4.5}
    \Phi_{1}(r_{\infty})= \lambda^{-(n-1)}(r_{\infty}).
\end{align}
Substituting (\ref{4.5}) into (\ref{4.3}), we get
\begin{align*}
    \int_{\Sigma_{\infty}}f^{\frac{n}{n-1}}d\mu_{\infty}= \omega_{n}.
\end{align*}
Hence,
\begin{align*}
    \int_{M_{0}} f^{\frac{n}{n-1}} d\mu \geq   \int_{\Sigma_{0}} f^{\frac{n}{n-1}} d\mu \geq  \int_{\Sigma_{\infty}}f^{\frac{n}{n-1}}d\mu_{\infty} \geq \omega_{n},
\end{align*}
which implies that the starshaped hypersurface satisfies the inequalities (\ref{1.10}) and (\ref{1.11}).

It is obvious that the equality holds in the inequality (\ref{1.11}) for geodesic spheres, we just need to prove the converse. Suppose that the smooth starshaped hypersurface $\Sigma_{t}$ makes the equality hold, i.e.,
\begin{align}\label{4.6}
    \int_{\Sigma_{t}}f^{\frac{n}{n-1}}d\mu_{t}= \Phi_{1}^{\frac{n}{n-1}}(r_{\infty})\lambda^{n}(r_{\infty})\omega_{n}=\omega_{n}.
\end{align}
Then, along the flow (\ref{1.9}), we have $\frac{\partial}{\partial t}\int_{\Sigma_{t}}f^{\frac{n}{n-1}}d\mu_{t} =0$. From (\ref{4.2}), we get 
\begin{align*}
  \frac{n}{n-1}\frac{\partial \Phi_{1}}{\partial r}\sqrt{1+\lambda^{-2}|D r|^{2}} +\Phi_{1} H = 0.
\end{align*}
Thereby, 
\begin{align*}
    \frac{\partial r}{\partial t} = -\frac{n}{n-1}\frac{\partial \Phi_{1}}{\partial r}(1+\lambda^{-2}|D r|^{2})-\Phi_{1} H  =0,
\end{align*}
and
\begin{align*}
    \frac{\partial}{\partial t}\left(\Phi_{1}^{\frac{n}{n-1}}\lambda^{n}\right) = \Phi_{1}^{\frac{1}{n-1}}\lambda^{n-1}\left(\frac{n}{n-1}\frac{\partial \Phi_{1}}{\partial r}\lambda+n\Phi_{1}\lambda^{'}\right)\frac{\partial r}{\partial t} =0.
\end{align*}
Thus,
\begin{align*}
    \lambda^{n}\Phi_{1}^{\frac{n}{n-1}}(r) =\lambda^{n}(r_{\infty})\Phi_{1}^{\frac{n}{n-1}}(r_{\infty})
\end{align*}
and the equation (\ref{4.6}) is equivalent to
\begin{align}\label{4.7}
    \int_{\Sigma_{t}}\Phi_{1}^{\frac{n}{n-1}}d\mu_{t}= \Phi_{1}^{\frac{n}{n-1}}(r_{\infty})\lambda^{n}(r_{\infty})\omega_{n}=\Phi_{1}^{\frac{n}{n-1}}(r)\lambda^{n}(r)\omega_{n}.
\end{align}
By (\ref{4.7}), it follows that
\begin{align}\label{4.8}
    \lambda^{n}(r) = \frac{\int_{\Sigma_{t}}\Phi_{1}^{\frac{n}{n-1}}d\mu_{t}}{\Phi_{1}^{\frac{n}{n-1}}(r)\omega_{n}}.
\end{align}
Differentiating (\ref{4.8}) and combining with (\ref{4.7}), we get
\begin{align*}
    \nabla (\lambda^{n}) = \nabla \left( \frac{\int_{M_{t}}\Phi_{1}^{\frac{n}{n-1}}d\mu_{t}}{\Phi_{1}^{\frac{n}{n-1}}(r)\omega_{n}} \right)
    = - \frac{n}{n-1}\frac{\nabla \Phi_{1}}{\Phi_{1}}\lambda^{n},
\end{align*}
i.e.,
\begin{align*}
    \nabla(\lambda^{n})+\frac{n}{n-1}\frac{\nabla \Phi_{1}}{\Phi_{1}}\lambda^{n}
    = \nabla r \left(\frac{n}{n-1}\frac{\partial \Phi_{1}}{\partial r}+n\Phi_{1}\frac{\lambda^{'}}{\lambda}\right)\lambda^{n}\Phi_{1}^{-1} = 0.
\end{align*}
It can be inferred that either $\nabla r = 0$ or $\frac{n}{n-1}\frac{\partial \Phi_{1}}{\partial r}+n\Phi_{1}\frac{\lambda^{'}}{\lambda}=0$.

Case 1: If $\nabla r = 0 $, $r$ is constant, then $\Sigma_{t}$ is a geodesic sphere.

Case 2: If $\frac{n}{n-1}\frac{\partial \Phi_{1}}{\partial r}+n\Phi_{1}\frac{\lambda^{'}}{\lambda}=0$, then
\begin{align}\label{4.9}
    \Phi_{1}(r)=\lambda^{-(n-1)}(r).
\end{align}
Substituting (\ref{4.9}) into (\ref{4.4}), we have
\begin{align*}
    \lambda^{-n}\left(\lambda H -n\lambda^{'}v \right)= 0.
\end{align*}
Then
\begin{align}\label{4.10}
    H= \frac{n\lambda^{'}v}{\lambda}.
\end{align}
Let $\alpha^{ki}=\frac{1}{\lambda v}\left(\sigma^{ki}-\frac{\varphi^{k}\varphi^{i}}{v^{2}}\right)$. By (\ref{2.5})and (\ref{4.10}), we deduce that $\alpha^{ki}\varphi_{ik}=\frac{n\lambda^{'}v}{\lambda}-\frac{n\lambda^{'}}{\lambda v}$. Also, from (\ref{2.4}), it follows that
\begin{align*}
    h^{i}_{j} = \frac{\lambda^{'}}{\lambda v}\delta^{i}_{j}-\frac{1}{\lambda v}\left(\sigma^{ki}-\frac{\varphi^{k}\varphi^{i}}{v^{2}}\right)\varphi_{kj}
    =\frac{\lambda^{'}}{\lambda v}\delta^{i}_{j}-\alpha^{ki}\varphi_{kj}
    =\frac{\lambda^{'}}{\lambda v}\delta^{i}_{j}-\left(\frac{\lambda^{'}}{\lambda v}-\frac{\lambda^{'} v}{\lambda }\right)\delta^{i}_{j}.
\end{align*}
Thus,
\begin{align*}
    |A|^{2}=h^{i}_{j}h^{j}_{i} = n\left(\frac{\lambda^{'} v}{\lambda}\right)^{2}.
\end{align*}
Combining (\ref{4.10}) and the above equality, we have 
\begin{align*}
    \frac{|A|^{2}}{H^{2}} = \frac{1}{n}.
\end{align*}
Therefore, $\Sigma_{t}$ is a geodesic sphere centered at the origin for $t>0$. Since $\Sigma_{0}$ can be approximated smoothly by $\Sigma_{t}$, i.e., $\Sigma_{0}=\lim_{t \to 0} \Sigma_{t} $, $\Sigma_{0}$ is also  a geodesic sphere centered at the origin. This completes the proof of Theorem \ref{th1.12}.
\hfill${\square}$

\subsection{ Sharp Michael-Simon type inequality for
$k$-th mean curvatures  }
\
\vglue-10pt
 \indent

In this section, we will give the proof of the new geometric inequalities (\ref{1.14}) and (\ref{1.15}) by applying  the convergence result of the flow (\ref{1.12}).

\noindent {\bf Proof of Theorem \ref{th1.15}}: First simplify the inequality (\ref{1.14}). We may assume that
\begin{align*}
    \int_{M_{0}} \lambda^{'} \sqrt{f^{2}E_{k}^{2} +|\nabla^{M_{0}} f|^{2} E_{k-1}^{2}} -\int_{M_{0}}\left\langle \bar{\nabla}\left(f\lambda^{'}\right),\nu \right\rangle  E_{k-1} &+\int_{\partial M_{0}}  f  E_{k-1}  \\
    &=\int_{M_{0}}f^{\frac{n-k+1}{n-k}} E_{k-1}.
\end{align*}
This normalisation ensures the existence of solution $\eta: M_{0} \rightarrow \mathbb{R}$ for the following PDE
\begin{align*}
    \operatorname{div}_{M_{0}}\left(E_{k-1}f \nabla^{M} \eta \right)=f^{\frac{n-k+1}{n-k}} E_{k-1}- \lambda^{'} \sqrt{f^{2}E_{k}^{2} +|\nabla^{M_{0}} f|^{2} E_{k-1}^{2}} + \left\langle \bar{\nabla}\left(f\lambda^{'}\right),\nu \right\rangle  E_{k-1}
\end{align*}
on $M_{0}$, and $\left\langle \nabla^{M_{0}} \eta , \vec{n} \right\rangle=1$ on $\partial M_{0}$. Here, $\vec{n}$ denotes the co-normal to $M_{0}$.  According to the standard elliptic regularity theory, $\eta \in C^{2, \beta}$ for each $0<\beta<1$. Hence, (\ref{1.14}) can be reduced into the following inequality
\begin{align}\label{4.11}
    \int_{M_{0}}E_{k-1}f^{\frac{n-k+1}{n-k}} d\mu \geq \int_{\Sigma_{0}}E_{k-1}f^{\frac{n-k+1}{n-k}} d\mu \geq p_{k}\circ h_{0}^{-1}\left(W^{\lambda^{'}}_{0}(\Omega_{0})\right).
\end{align}

Secondly, we prove that $W^{\lambda^{'}}_{0}(\Omega_{t})$ monotonically increases along the flow (\ref{1.12}), where $\Omega_{t}$ is the domain enclosed by $\Sigma_{t}$. Using (\ref{2.10}) and (\ref{2.20}), we have
\begin{align*}
    \frac{\partial}{\partial t}W^{\lambda^{'}}_{0}(\Omega_{t})=(n+1)\int_{\Sigma_{t}}\left(\lambda^{'}\frac{E_{k-2}}{E_{k-1}}-u\right)d\mu_{t}
    \geq (n+1)\int_{\Sigma_{t}}\left(\frac{\lambda^{'}}{E_{1}}-u\right)d\mu_{t} \geq 0,
\end{align*}
where the Heintze-Karcher inequality (see Theorem 3.5 in \cite{B13}) is used in the last inequality.

Finally we prove the monotonicity of $\int_{\Sigma_{t}}E_{k-1}f^{\frac{n-k+1}{n-k}} d\mu_{t}$ along the flow (\ref{1.12}). For the  convenience of computation, we use $\widehat{\Phi}_{2}(r)=\Phi_{2}^{\frac{n-k+1}{n-k}}(r(\xi,t))$. Along the general flow (\ref{2.11}), we have
\begin{align*}
    \frac{\partial}{\partial t} \int_{\Sigma_{t}}  E_{k-1} \widehat{\Phi}_{2} d\mu_{t}= &\int_{\Sigma_{t}}  \widehat{\Phi}_{2} \dot{E}^{ij}_{k-1}\left(-\nabla_{j}\nabla_{i}F - F(h^{2})_{ij}+Fg_{ij}\right)d\mu_{t} \\
    &+ \int_{\Sigma_{t}} \frac{\partial \widehat{\Phi}_{2}}{\partial r}\frac{\partial r}{\partial t} E_{k-1} d\mu_{t}
    + \int_{\Sigma_{t}} nE_{k-1}E_{1}\widehat{\Phi}_{2} F d\mu_{t}.
\end{align*}
Since $\dot{E}^{ij}_{k-1}$ is divergence-free, we have 
\begin{align*}
    \frac{\partial}{\partial t} \int_{\Sigma_{t}}  E_{k-1} \widehat{\Phi}_{2} d\mu_{t}=
    &\int_{\Sigma_{t}}-\frac{k-1}{n}\Delta_{\Sigma_{t}}\widehat{\Phi}_{2} E_{k-2} F d\mu_{t} + \int_{\Sigma_{t}} (k-1)\widehat{\Phi}_{2}E_{k-2}F d\mu_{t}\\
    &+ \int_{\Sigma_{t}} \frac{\partial \widehat{\Phi}_{2}}{\partial r} v E_{k-1}F d\mu_{t}+ \int_{\Sigma_{t}} (n-k+1)\widehat{\Phi}_{2}E_{k}F d\mu_{t}\\
    =& \int_{\Sigma_{t}} \left(-\frac{k-1}{n}\Delta_{\Sigma_{t}} \widehat{\Phi}_{2}+(k-1)\widehat{\Phi}_{2}\right) E_{k-2} F d\mu_{t} \\
    &+ \int_{\Sigma_{t}} \frac{\partial \widehat{\Phi}_{2}}{\partial r} v E_{k-1}F d\mu_{t}+ \int_{\Sigma_{t}} (n-k+1)\widehat{\Phi}_{2}E_{k}F d\mu_{t},
\end{align*}
where we use integration by parts, (\ref{2.6}) and (\ref{2.8}). From Assumption \ref{as1.14} (1), it follows that  
\begin{align*}
    -\frac{k-1}{n}\Delta_{\Sigma_{t}} \widehat{\Phi}_{2}+(k-1)\widehat{\Phi}_{2} =-\frac{\partial \widehat{\Phi}_{2}}{\partial r} \frac{\lambda^{'}}{\lambda}v^{2} + m\frac{E_{k-1}}{E_{k-2}} \lambda (\lambda^{'})^{m-k}.
\end{align*} 
Thus
\begin{align*}
    \frac{\partial}{\partial t} &\int_{\Sigma_{t}}  E_{k-1} \widehat{\Phi}_{2} d\mu_{t} \\
    =&\int_{\Sigma_{t}}  m\lambda (\lambda^{'})^{m-k} E_{k-1} F d\mu_{t} - \int_{\Sigma_{t}}  \left(\frac{\partial \widehat{\Phi}_{2}}{\partial r}\frac{\lambda^{'}}{\lambda}v^{2}E_{k-2}
    - \frac{\partial \widehat{\Phi}_{2}}{\partial r} v E_{k-1}\right)F d\mu_{t} \\
    &+ \int_{\Sigma_{t}} (n-k+1)\widehat{\Phi}_{2}E_{k}F d\mu_{t} .
\end{align*}
Substituting $F =\frac{E_{k-2}}{E_{k-1}}-\frac{u}{\lambda^{'}}$ into the above equation and applying the Newton-MacLaurin inequality (\ref{2.10}), we can obtain the following along the flow (\ref{1.12})
\begin{align}\label{4.12}
    \frac{\partial}{\partial t} \int_{\Sigma_{t}}  E_{k-1} \widehat{\Phi}_{2} d\mu_{t} 
    =&\int_{\Sigma_{t}}  m\lambda (\lambda^{'})^{m-k-1} \left(\lambda^{'}E_{k-2} - uE_{k-1} \right)  - \frac{\partial \widehat{\Phi}_{2}}{\partial r}\frac{\lambda^{'}}{\lambda}v^{2}E_{k-1}\left(\frac{E_{k-2}}{E_{k-1}}-\frac{u}{\lambda^{'}}\right)^{2} d\mu_{t} \notag \\
    &+ \int_{\Sigma_{t}} (n-k+1)\widehat{\Phi}_{2}E_{k}\left(\frac{E_{k-2}}{E_{k-1}}-\frac{u}{\lambda^{'}}\right) d\mu_{t}  \notag \\
    \leq& \int_{\Sigma_{t}}  m\lambda (\lambda^{'})^{m-k-1}  \left(\lambda^{'}E_{k-2} - uE_{k-1} \right)  -   \frac{\partial \widehat{\Phi}_{2}}{\partial r}\frac{\lambda^{'}}{\lambda}v^{2}E_{k-1}\left(\frac{E_{k-2}}{E_{k-1}}-\frac{u}{\lambda^{'}}\right)^{2} d\mu_{t}
    \\
    &+\int_{\Sigma_{t}} (n-k+1)\frac{\widehat{\Phi}_{2}}{\lambda^{'}}\left(\lambda^{'}E_{k-1} - uE_{k}\right) d\mu_{t}. \notag
\end{align}
From Assumption \ref{as1.14}(2), we have 
\begin{align*}
    \frac{\partial \widehat{\Phi}_{2}}{\partial \lambda^{'}} - \frac{\widehat{\Phi}_{2}}{\lambda^{'} }\geq 0,\qquad
    \frac{\partial \widehat{\Phi}_{2}}{\partial r} = \frac{\partial \widehat{\Phi}_{2}}{\partial \lambda^{'}}  \lambda \geq \widehat{\Phi}_{2}\frac{\lambda}{\lambda^{'} } >0,
\end{align*} 
and combining (\ref{2.2}), (\ref{2.6}) and (\ref{2.7}), it follows that 
\begin{align*}
    \lambda^{'}E_{k-1}-uE_{k}=\frac{1}{k}\dot{E}^{ij}_{k}\left(\lambda^{'}g_{ij}-uh_{ij}\right)=\frac{1}{k}\dot{E}^{ij}_{k}\nabla_{i}\nabla_{j}\lambda^{'},
\end{align*}
Thus
\begin{align*}  
    \frac{\partial}{\partial t} &\int_{\Sigma_{t}}  E_{k-1} \widehat{\Phi}_{2} d\mu_{t} \\ \leq&  \int_{\Sigma_{t}} -\frac{m}{k-1} \dot{E}^{ij}_{k-1}\nabla_{i}(\lambda (\lambda^{'})^{m-k-1}) \nabla_{j}\lambda^{'}   -  \frac{n-k+1}{k}  \dot{E}^{ij}_{k}\nabla_{i}\frac{\widehat{\Phi}_{2}}{\lambda^{'}}\nabla_{j}\lambda^{'} d\mu_{t}
    \\
    =& \int_{\Sigma_{t}} -\frac{m}{k-1} (\lambda^{'})^{m-k-2}\lambda\left[\frac{(\lambda^{'})^{2}}{\lambda^{2}}+ (m-k-1)\right]\dot{E}^{ij}_{k-1}\nabla_{i}\lambda^{'} \nabla_{j}\lambda^{'} \\
    &- \frac{n-k+1}{k}\left(\frac{\partial \widehat{\Phi}_{2}}{\partial \lambda^{'}} \frac{1}{\lambda^{'}}-\frac{\widehat{\Phi}_{2}}{(\lambda^{'})^{2}}\right)\dot{E}^{ij}_{k}\nabla_{i}\lambda^{'}\nabla_{j}\lambda^{'} d\mu_{t} \\
    \leq& \int_{\Sigma_{t}} -\frac{m(m-k)}{k-1} (\lambda^{'})^{m-k-2}\lambda \dot{E}^{ij}_{k-1}\nabla_{i}\lambda^{'} \nabla_{j}\lambda^{'} \\
    &- \frac{n-k+1}{k}\left(\frac{\partial \widehat{\Phi}_{2}}{\partial \lambda^{'}} \frac{1}{\lambda^{'}}-\frac{\widehat{\Phi}_{2}}{(\lambda^{'})^{2}}\right)\dot{E}^{ij}_{k}\nabla_{i}\lambda^{'}\nabla_{j}\lambda^{'} d\mu_{t}\\
    \leq&  0 .
\end{align*}
where the last inequality follows from the positivity of $\dot{E}^{ij}_{k}$ and $\dot{E}^{ij}_{k-1}$. Therefore
\begin{align*}
    \int_{M_{0}} E_{k-1}\widehat{\Phi}_{2} d\mu \geq \int_{\Sigma_{0}} E_{k-1}\widehat{\Phi}_{2} d\mu_{t} \geq \int_{\Sigma_{\infty}}E_{k-1}\widehat{\Phi}_{2} d\mu_{\infty}=\int_{B_{R}}E_{k-1}\widehat{\Phi}_{2} d\mu_{\mathbb{H}^{n+1}},
\end{align*}
where the last equality is obtained from the convergence result of the flow (\ref{1.12}) and $B_{R}=\partial B^{n+1}_{R}$. Also $\nabla \Phi_{2}=\frac{\partial \Phi_{2}}{\partial r} \nabla r=0$ on $B_{R}$, i.e., $f = \Phi_{2}(r)$ is constant on $B_{R}$. Hence
\begin{align*}
    \int_{B_{R}}  E_{k-1} \widehat{\Phi}_{2} d\mu_{\mathbb{H}^{n+1}} &=\omega_{n} \Phi_{2}^{\frac{n-k+1}{n-k}}\int_{B_{R}} E_{k-1}(\kappa) d\mu_{\mathbb{H}^{n+1}}= \omega_{n} \Phi_{2}^{\frac{n-k+1}{n-k}}(\lambda^{'})^{k-1}\lambda^{n-k+1} (R)\\
    &=p_{k}\circ h^{-1}_{0}\left(W^{\lambda^{'}}_{0}(B^{n+1}_{R})\right). \notag
\end{align*}
Furthermore, we already know that $W^{\lambda^{'}}_{0}(\Omega_{t})$ is monotonically increasing along the flow (\ref{1.12}) and $p_{k}(r)$ is monotonically increasing with respect to $r$, thus 
\begin{align*}
    \int_{M_{0}} E_{k-1} \widehat{\Phi}_{2} d\mu_{t} &\geq \int_{\Sigma_{0}} E_{k-1} \widehat{\Phi}_{2} d\mu_{t} \\
    &\geq\int_{B_{R}}  E_{k-1} \widehat{\Phi}_{2} d\mu_{\mathbb{H}^{n+1}} = p_{k}\circ h^{-1}_{0}\left(W^{\lambda^{'}}_{0}(B^{n+1}_{R})\right)
    \geq p_{k}\circ h^{-1}_{0}\left(W^{\lambda^{'}}_{0}(\Omega_{0})\right),
\end{align*}
which means that the starshaped, strictly $k$-convex hypersurface satisfies the
inequalities (\ref{1.14}) and (\ref{1.15}).

We now prove that $\Sigma_{0}$ is a geodesic sphere when the equality holds in (\ref{1.15}). If a smooth starshaped, strictly $k$-convex hypersurface $\Sigma_{t}$ attains the equality
\begin{align*}
    \int_{\Sigma_{t}}E_{k-1}f^{\frac{n-k+1}{n-k}} d\mu_{t} = p_{k}\circ h^{-1}_{0}\left(W^{\lambda^{'}}_{0}(\Omega_{t})\right).
\end{align*}
Then, $\frac{\partial}{\partial t}\int_{\Sigma_{t}}E_{k-1}f^{\frac{n-k+1}{n-k}} d\mu_{t}=0$. According to Lemma \ref{lem2.4} and (\ref{4.12}), $\Sigma_{t}$ is a geodesic sphere centered at the origin for $t>0$. Since $\Sigma_{0}$ can be smoothly approximated by a family of geodesic spheres, $\Sigma_{0}$ is also a geodesic sphere. This completes the proof of Theorem \ref{th1.15}.

\hfill${\square}$



\end{document}